\documentclass{birkart}
\usepackage{amssymb}

\expandafter\ifx\csname delta.sty\endcsname\relax \else\endinput\fi
\expandafter\edef\csname delta.sty\endcsname{%
 \catcode`\noexpand\@=\the\catcode`\@\space}

\let\atbefore @

\catcode`\@=11

\newif\ifMag
\let\@ft@\expandafter
\numberwithin{equation}{section}

\newif\ifComments

\def\tod@y{\ifcase\month\or
 January\or February\or March\or April\or May\or June\or July\or
 August\or September\or October\or November\or December\fi\space\,
\number\day,\space\,\number\year}

\def\h@@r{hh}\def\m@n@te{mm}
\def\wh@tt@me{\count@\time\divide\count@ 60\edef\h@@r{\number\count@}%
 \multiply\count@ -60\advance\count@\time\edef
 \m@n@te{\ifnum\count@<10 0\fi\number\count@}}
\def\t@me{\h@@r\/{\rm:}\m@n@te} \let\whattime\wh@tt@me
\let\Today\tod@y \let\nowtime\t@me

\def\ftext#1{{\let\thefootnote\relax\footnotetext{\vsk-.8>\nt #1}}}
\def\em#1{{\itshape #1\/}}

 {\end{trivlist}\egroup\global\@ignoretrue\ifComments\unvbox\commentbox\fi
 \noindent}

\def\gadv{\global\adv} \def\gad#1{\gadv#1\@ne} \def\gadneg#1{\gadv#1-\@ne}
\def\textindent#1{\indent\llap{#1\enspace}\ignorespaces}

\newcount\itemlet
\def\newbi{\itemlet 96} \newbi
\def\bitem{\gad\itemlet\endgraf\hangindent1.5\parindent
 \hglue-.5\parindent\textindent{\upshape\rlap{\char\the\itemlet}\hp{b})}}

\newcount\itemrm

\def\iitem{\gad\itemrm\endgraf\hangindent1.5\parindent\hglue-.5\parindent
 \textindent{\upshape\hp{v}\llap{\romannumeral\the\itemrm})}}

\let\Disp\[ \let\endD\] \let\{\protect

\def\Tag#1{\label{e:#1}\let\notag\relax} 
\def\sh@nd#1#2{\begin{#1*}#2\end{#1*}}
\def\n@t@gs#1#2#3{\let\n@@@l\\ \begin{#1}#2\global\let\d@bl\\
 \gdef\\{\notag\d@bl}#3\notag\global\let\\\d@bl\end{#1}\let\\\n@@@l}

\def\Gather#1\endG{\n@t@gs{gather}{}{#1}}
\def\gAther#1\endG{\sh@nd{gather}{#1}}
\def\Align#1\endA{\n@t@gs{align}{}{#1}}
\def\aLign#1\endA{\sh@nd{align}{#1}}
\def\Alignat#1#2\endAt{\n@t@gs{alignat}{#1}{#2}}
\def\aLignat#1\endAt{\sh@nd{alignat}{#1}}

\def\(#1){\textup{(\ref{e:#1})}}
\def\[{\@ifnextchar:\c@t@sect\c@t@}
\def\c@t@sect:#1]{\ref{s:#1}} \def\c@t@#1]{\ref{t:#1}}

\def\qed{\hbox{}\nobreak\hfill\nobreak{\m@th$\,\square$}}
\def\sk@@p#1{\par\skip@#1\relax\ifdim\lastskip<\skip@\relax\removelastskip
 \vskip\skip@\fi}
\def\demo#1{\sk@@p\medskipamount\nt{\ignore\it #1\unskip.}\enspace
 \ignore}
\def\enddemo{\sk@@p\medskipamount}
 
\def\Pf#1.{\demo{Proof #1}} 

\def\Text#1{\crcr\noalign{\alb\vsk>\normalbaselines\vsk->\vbox{\nt #1\strut}%
 \nobreak\nointerlineskip\vbox{\strut}\nobreak\vsk->\nobreak}}

 \let\nl\newline
\let\bls\baselineskip \let\ignore\ignorespaces \let\adv\advance
\def\vsk#1>{\vskip#1\bls}
\def\vv#1>{\vadjust{\vsk#1>}\ignore}
\def\vvn#1>{\vadjust{\nobreak\vsk#1>\nobreak}\ignore}
\def\vvv#1>{\vskip\z@\vsk#1>\nt\ignore}
\def\vvgood{\vadjust{\penalty-500}} 

\def\mathbox#1{\hbox{\m@th$#1$}}

\let\dsize\displaystyle \let\tsize\textstyle
\let\ssize\scriptstyle \let\sss\scriptscriptstyle
 
\let\vp\vphantom \let\hp\hphantom \let\nt\noindent
\def\hline{\hbox to\hsize}
\let\cline\centerline \let\lline\leftline \let\rline\rightline
\def\nn#1>{\noalign{\vskip#1\p@}} \def\NN#1>{\openup#1\p@}
 
\let\Lim\lim \def\lim{\Lim\limits} \let\Sum\sum \def\sum{\Sum\limits}
\def\Plus{\bigoplus\limits} 
\let\Prod\prod \def\prod{\Prod\limits} \let\Int\int \def\int{\Int\limits}

\def\~{\leavevmode\@ifnextchar~\m@n@s\@md@sh}
\def\m@n@s~{\raise.15ex\mathbox{-}} \def\@md@sh{\raise.13ex\hbox{--}}
\let\procent\% \def\%#1{\ifmmode\mathop{#1}\limits\else\procent#1\fi}
\let\@ml@t\" \def\"#1{\ifmmode ^{(#1)}\else\@ml@t#1\fi}
\let\@c@t@\' \def\'#1{\ifmmode _{(#1)}\else\@c@t@#1\fi}
\let\colon\: \def\:{^{\vp|}} \def\&{.\kern.1em} \def\^#1{\text{\m@th#1}}

\newif\ifNewskips

\def\Newskips{\global\Newskipstrue
 \gdef\>{\relax\ifmmode\mskip.666667\thinmuskip\relax\else\kern.111111em\fi}
 \gdef\}{\relax\ifmmode\mskip-.666667\thinmuskip\relax\else\kern-.111111em\fi}
 \gdef\){\relax\ifmmode\mskip.333333\thinmuskip\relax\else\kern.0555556em\fi}
 \gdef\]{\relax\ifmmode\mskip-.333333\thinmuskip\relax\else\kern-.0555556em\fi}}
\Newskips

\def\Re{\mathop{\mathrm{Re}\>}} \def\Im{\mathop{\mathrm{Im}\>}}
\def\End{\mathop{\mathrm{End}\>}} 
 \def\Sym{\mathop{\mathrm{Sym}\)}}

\def\id{\mathrm{id}}  
\def\1{^{-1}} \def\_#1{_{\Rlap{#1}}}
\def\vst#1{{\lower1.9\p@
 \hbox{\m@th$\bigr|_{\raise.5\p@\hbox{\m@th$\ssize#1$}}$}}}
\def\vrp#1:#2>{{\vrule height#1 depth#2 width\z@}}
\def\vru#1>{\vrp#1:\z@>} \def\vrd#1>{\vrp\z@:#1>}
\def\qqq{\qquad\quad} 
\def\sscr#1{\raise.3ex\hbox{\m@th$\sss#1$}} \def\@@PS{\mathbf{OOPS!!!}}

\def\lsym#1{#1\alb\ldots\relax#1\alb}
\def\lc{\lsym,}  \def\lx{\lsym\x} \def\lox{\lsym\ox}

\let\texspace\ \def\ {\ifmmode\alb\fi\texspace}

\def\Line#1{\kern-.5\hsize\hline{\m@th$\dsize#1$}\kern-.5\hsize}
\def\Lline#1{\kern-.5\hsize\lline{\m@th$\dsize#1$}\kern-.5\hsize}
\def\Cline#1{\kern-.5\hsize\cline{\m@th$\dsize#1$}\kern-.5\hsize}
\def\Rline#1{\kern-.5\hsize\rline{\m@th$\dsize#1$}\kern-.5\hsize}

\def\Ll@p#1{\llap{\m@th$#1$}} \def\Rl@p#1{\rlap{\m@th$#1$}}
 \def\Cl@p#1{\llap{\m@th$#1$\hss}}
\def\Llap#1{\mathchoice{\Ll@p{\dsize#1}}{\Ll@p{\tsize#1}}{\Ll@p{\ssize#1}}%
 {\Ll@p{\sss#1}}}
\def\Clap#1{\mathchoice{\Cl@p{\dsize#1}}{\Cl@p{\tsize#1}}{\Cl@p{\ssize#1}}%
 {\Cl@p{\sss#1}}}
\def\Rlap#1{\mathchoice{\Rl@p{\dsize#1}}{\Rl@p{\tsize#1}}{\Rl@p{\ssize#1}}%
 {\Rl@p{\sss#1}}}
 
\def\LRtph#1#2{\setbox\z@\hbox{#1}\dimen\z@\wd\z@\hbox{\hbox to\dimen\z@{#2}}}
\def\LRph#1#2{\LRtph{\m@th$#1$}{\m@th$#2$}}

\def\LLph#1#2{\LRph{#1}{\hss#2}} 

\def\Lph#1#2{\mathchoice{\LLph{\dsize#1}{\dsize#2}}{\LLph{\tsize#1}{\tsize#2}}
 {\LLph{\ssize#1}{\ssize#2}}{\LLph{\sss#1}{\sss#2}}}

\def\Lto#1{\setbox\z@\hbox{\m@th$\tsize{#1}$}%
 \mathrel{\mathop{\hbox to\wd\z@{\rightarrowfill}}\limits#1}}
\def\Lgets#1{\setbox\z@\hbox{\m@th$\tsize{#1}$}%
 \mathrel{\mathop{\hbox to\wd\z@{\leftarrowfill}}\limits#1}}

\def\vpb#1{{\vp{\big(}}^{\]#1}}

\let\alb\allowbreak

\let\o\circ \let\x\times \let\ox\otimes 
\let\sub\subset 
\let\le\leqslant \let\ge\geqslant
\let\der\partial \let\8\infty \let\*\star
\let\bra\langle \let\ket\rangle
 
\let\map\mapsto  \let\hto\hookrightarrow
 
  \def\nin{\not\in}

\let\lb\lbrace \let\rb\rbrace

\let\Bbb\mathbb

\def\pms{\raise.25ex\mathbox{\ssize\pm}\>}
\def\mps{\raise.25ex\mathbox{\ssize\mp}\>}

\let\al\alpha
\let\bt\beta
\let\gm\gamma \let\Gm\Gamma 
\let\dl\delta  
 \let\eps\varepsilon \let\epsilon\eps

\let\zt\zeta

\let\ka\kappa
\let\la\lambda 

\let\si\sigma 
 
\let\pho\phi \let\phi\varphi

\let\om\omega \let\Om\Omega 

\def\C{\Bbb C}
\def\R{\Bbb R}
\def\Z{\Bbb Z}

\def\AA{\Bbb A}

\def\II{\Bbb I}

\def\Zp{\Z_{\ge 0}} \def\Zpp{\Z_{>0}}

\def\h@ph{\discretionary{}{}{-}} \def\$#1$-{\,\^{$#1$}\h@ph}

\def\difl/{differential} \def\dif/{difference}
\def\cf.{cf.\ \ignore} \def\Cf.{Cf.\ \ignore}
\def\egv/{eigenvector} \def\eva/{eigenvalue} \def\eq/{equation}
\def\lhs/{the left hand side} \def\rhs/{the right hand side}
\def\Lhs/{The left hand side} \def\Rhs/{The right hand side}
\def\gby/{generated by} \def\wrt/{with respect to} \def\st/{such that}
\def\resp/{respectively} \def\off/{offdiagonal} \def\wt/{weight}
\def\pol/{polynomial} \def\rat/{rational} \def\tri/{trigonometric}
\def\fn/{function} \def\var/{variable} \def\raf/{\rat/ \fn/}
\def\inv/{invariant} \def\hol/{holomorphic} \def\hof/{\hol/ \fn/}
\def\mer/{meromorphic} \def\mef/{\mer/ \fn/} \def\mult/{multiplicity}
\def\sym/{symmetric} \def\perm/{permutation}
\def\rep/{representation} \def\irr/{irreducible} \def\irrep/{\irr/ \rep/}
\def\hom/{homomorphism} \def\aut/{automorphism} \def\iso/{isomorphism}
\def\lex/{lexicographical} \def\as/{asymptotic} \def\asex/{\as/ expansion}
\def\ndeg/{nondegenerate} \def\neib/{neighbourhood} \def\deq/{\dif/ \eq/}
\def\hw/{highest \wt/} \def\gv/{generating vector} \def\eqv/{equivalent}
\def\msd/{method of steepest descend} \def\pd/{pairwise distinct}
\def\wlg/{without loss of generality} \def\Wlg/{Without loss of generality}
\def\onedim/{one-dim\-en\-sion\-al} \def\fd/{fi\-ni\-te-dim\-en\-sion\-al}
\def\qcl/{quasiclassical} \def\hwv/{\hw/ vector}
\def\hgeom/{hyper\-geo\-met\-ric} \def\hint/{\hgeom/ integral}
\def\hwm/{\hw/ module} \def\emod/{evaluation module} \def\Vmod/{Verma module}
\def\symg/{\sym/ group} \def\sol/{solution} \def\eval/{evaluation}
\def\anf/{analytic \fn/} \def\anco/{analytic continuation}
\def\qg/{quantum group} \def\qaff/{quantum affine algebra}

\def\Rm/{\^{$R$-}matrix} \def\Rms/{\^{$R$-}matrices}
\def\YB/{Yang-Baxter \eq/}
\def\Ba/{Bethe ansatz} \def\Bv/{Bethe vector} \def\Bae/{\Ba/ \eq/}
\def\KZv/{Knizh\-nik-Zamo\-lod\-chi\-kov} \def\KZvB/{\KZv/-Bernard}
\def\KZ/{{\sl KZ\/}} \def\qKZ/{{\sl qKZ\/}}
\def\KZB/{{\sl KZB\/}} \def\qKZB/{{\sl qKZB\/}}
\def\qKZo/{\qKZ/ operator} \def\qKZc/{\qKZ/ connection}
\def\KZe/{\KZ/ \eq/} \def\qKZe/{\qKZ/ \eq/} \def\qKZBe/{\qKZB/ \eq/}

\def\LPT/{Laboratoire de Physique Th\'eorique ENSLAPP}
\def\ENSLyon/{\'Ecole Normale Sup\'erieure de Lyon}
\def\DMS/{Department of Mathematics, Faculty of Science}
\def\DMO/{\DMS/, Osaka University}
\def\DMOaddr/{Toyonaka, Osaka 560, Japan}
\def\dmoemail/{vt@math.sci.osaka-u.ac.jp}
\def\SPb/{St\&Petersburg}
\def\home/{\SPb/ Branch of Steklov Mathematical Institute}
\def\homeaddr/{Fontanka 27, \SPb/ \,191011, Russia}
\def\homemail/{vt@pdmi.ras.ru}
\def\absence/{On leave of absence from \home/}
\def\UNC/{Department of Mathematics, University of North Carolina}
\def\ChH/{Chapel Hill}
\def\UNCaddr/{\ChH/, NC 27599, USA} \def\avemail/{av@math.unc.edu}
\def\grant/{NSF grant DMS\~9501290}	
\def\Grant/{Supported in part by \grant/}

\def\Aomoto/{K\&Aomoto}
\def\Dri/{V\]\&G\&Drin\-feld}
\def\Fadd/{L\&D\&Fad\-deev}
\def\Feld/{G\&Felder}
\def\Fre/{I\&B\&Fren\-kel}
\def\Gustaf/{R\&A\&Gustafson}
\def\Kazh/{D\&Kazhdan} \def\Kir/{A\&N\&Kiril\-lov}
\def\Kor/{V\]\&E\&Kore\-pin}
\def\Lusz/{G\&Lusztig}
\def\MN/{M\&Naza\-rov}
\def\Resh/{N\&Reshe\-ti\-khin} \def\Reshy/{N\&\]Yu\&Reshe\-ti\-khin}
\def\Skl/{E\&K\&Sklya\-nin}
\def\SchV/{V\]\&\]V\]\&Schecht\-man} \def\Sch/{V\]\&Schecht\-man}
\def\Takh/{L\&A\&Takh\-tajan}
\def\VT/{V\]\&Ta\-ra\-sov} \def\VoT/{V\]\&O\&Ta\-ra\-sov}
\def\Varch/{A\&\]Var\-chenko} \def\Varn/{A\&N\&\]Var\-chenko}

\def\AMS/{Amer.\ Math.\ Society}
\def\CMP/{Comm.\ Math.\ Phys.{}}
\def\DMJ/{Duke.\ Math.\ J.{}}
\def\Inv/{Invent.\ Math.{}} 
\def\IMRN/{Int.\ Math.\ Res.\ Notices}
\def\JPA/{J.\ Phys.\ A{}}
\def\JSM/{J.\ Soviet\ Math.{}}
\def\LMP/{Lett.\ Math.\ Phys.{}}
\def\LMJ/{Leningrad Math.\ J.{}}
\def\LpMJ/{\SPb/ Math.\ J.{}}
\def\SIAM/{SIAM J.\ Math.\ Anal.{}}
\def\SMNS/{Selecta Math., New Series}
\def\TMP/{Theor.\ Math.\ Phys.{}}
\def\ZNS/{Zap.\ nauch.\ semin. LOMI}

\def\ASMP/{Advanced Series in Math.\ Phys.{}}

\def\AMSa/{AMS \publaddr Providence}
\def\Birk/{Birkh\"auser}
\def\CUP/{Cambridge University Press} \def\CUPa/{\CUP/ \publaddr Cambridge}
\def\Spri/{Springer-Verlag} \def\Spria/{\Spri/ \publaddr Berlin}
\def\WS/{World Scientific} \def\WSa/{\WS/ \publaddr Singapore}

\def\dash{\raise.13ex\hbox{-}}

\def\cnn#1>{\\[-\bls]\noalign{\vsk#1>}\notag}
\def\mmgood#1:#2>{\\\noalign{\vsk#1>\penalty-500\vsk#2>}\notag}

\def\fratop{\genfrac{}{}{0pt}1}
\def\q-{\>\hbox{\m@th$q\)$-}}
\def\+#1{^{\sscr{\bra\]#1\]\ket}}}
\def\##1{^{[#1]}}

\def\II{\mathbb I}
\def\PP{\mathbb P}
\def\VV{\mathbb V}

\def\db{\mathbf d}
\def\pb{\mathbf p}

\def\lbi{\^{\boldmath$l$}}
\def\mbi{\^{\boldmath$m$}}

\def\Cc{\mathcal C}

\def\Pc{\mathcal P}
\def\Zc{\mathcal Z}

\def\gg{\mathfrak g}
\def\hg{\mathfrak h}
\def\ng{\mathfrak n}
\def\gsl{\mathfrak{sl}}
\def\gl{\mathfrak{gl}}

\def\Ah{\widehat A}
\def\Dh{\widehat D}
\def\habla{\widehat\nabla}
\def\Phh{\hat\Phi}
\def\Uh{\widehat U}
\def\wh{\hat w}
\def\Xh{\widehat\Xi}

\def\At{\widetilde A}
\def\Ut{\widetilde U}
\def\Xt{\widetilde\Xi}

\def\zn{z_1\lc z_n} \def\lak{\la_1\lc\la_k}
\def\zlt{z_1,z_2;\alb\la_1,\la_2} \def\lzt{\la_1,\la_2;\alb z_1,z_2}
\def\xkn{x_{11}\lc x_{k1}\lc x_{1n}\)\lc x_{kn}}

\def\Vn{V_1\lc V_n} \def\Vox{V_1\]\lox V_n} \def\vox{v_1\]\lox v_n}

\def\ev{\mathop{\text{\sl ev}}\nolimits}
\def\lmb{\lbi\),\]\mbi} \def\mlb{\mbi\),\lbi}

\def\ngp{\ng_{\sss+}} \def\ngm{\ng_{\sss-}} 
\def\Pck{\Pc_k} \def\Pcm{\Pc_{\]m}} 
\def\Plm{\PP_{kn}[\)\lmb\)]} \def\Zlm{\Zc\)[\)\lmb\)]}
\def\Zknlm{\Zc_{kn}[\)\lmb\)]}
\def\zlla{\zn;\la_1,\la_2} \def\zllb{\zlla;\lbi} \def\zllm{\zllb\),\]\mbi}
\def\Udzlm{U_\db(\zn;\la_1,\la_2;\lmb)}
\def\Utdzlm{\Ut_\db(\zn;\la_1,\la_2;\lmb)}

\def\Ck{\C^{\)k}} \def\Cn{\C^{\)n}}
\def\glk{\gl_k} \def\gln{\gl_{\)n}} \def\glt{\gl_{\)2}}
\def\glkn{(\)\glk\>,\)\gln\))} 
\def\slk{\gsl_{\)k}} \def\sln{\gsl_{\)n}} \def\slt{\gsl_{\)2}}
\def\Uk{U(\glk)}  \def\Ukx{U\bigl(\glk[\)x\)]\bigr)}
\def\Yk{Y(\glk)} 
\def\Ukn{\bigl(U(\glk)\bigr)\vpb{\ox\)n}}

\def\Uqkn{\bigl(\)U_q(\glk)\>,\)U_q(\gln)\bigr)}

\def\Dlk{D_{\la_1}\lc D_{\la_k}}
\def\Dhk{\Dh_{\la_1}\lc\Dh_{\la_k}}
\def\Qlk{Q_{\la_1}\lc Q_{\la_k}}
\def\Zzn{Z_{z_1}\lc Z_{z_n}}
\def\nablasn{\nabla_{\}z_1}\lc\nabla_{\}z_n}}
\def\hablasn{\habla_{\}z_1}\lc\habla_{\}z_n}}
\def\Wla{W[\)\la\)]} 
\def\VVkn{\bigl(\)\VV\+k\)\bigr)\vpb{\ox n}}
\def\VVnk{\bigl(\)\VV\+n\)\bigr)\vpb{\ox k}}
\def\Vlm{(V_{l_1}\}\lox V_{l_n})\)[\)m_1,m_2]}
\def\Vlmt{(V_{l_1}\}\ox V_{l_2})\)[\)m_1,m_2]}
\def\Vklm{(V_{l_1}\+k\lox V_{l_n}\+k)\)[\)m_1\lc m_k]}
\def\Vnml{(V_{m_1}\+n\lox V_{m_k}\+n)\)[\)l_1\lc l_n]}
\def\Vmlt{(V_{m_1}\}\ox V_{m_2})\)[\)l_1, l_2]}

\def\Dmu{D_{\]\mu}}
\def\emu{e^{\>\mu}}

\def\Gmb#1{\Gm\bigl(#1\bigr)}
\def\Sin#1{\sin\)\bigl(#1\bigr)}
\def\sink#1{\Sin{\pi{(#1)}/\ka}}

\def\g-{\,\text{$\=\gg\)$-}}

\def\PVW{P_{VW}\:}
\def\RVV{R_{\)\VV\)\VV}\:} \def\RVlm{R_{V_lV_m\]}\:}
\def\RUV{R_{UV}\:} \def\RUW{R_{UW}\:} \def\RVW{R_{VW}\:}

\def\DD/{{\sl DD}}
\def\qDD/{{\sl q\)DD}}

\def\r-{\,\text{$\=r$-}}
\def\R-{\>\text{$\=R\)$-}}

\csname delta.sty\endcsname

\newtheorem{theorem}{Theorem}[section]
\newtheorem{proposition}[theorem]{Proposition}
\newtheorem{lemma}[theorem]{Lemma}

\newenvironment{remark}{\demo{\sl Remark}}{\enddemo}

\title[\smash{Duality for Knizhinik-Zamolodchikov Equations}]
{Duality for\\ Knizhinik-Zamolodchikov and\\ Dynamical Equations, and\\
Hypergeometric Integrals}

\author[\smash{V\].\,Tarasov}]{V\].\,Tarasov}

\address{St.\,Petersburg Branch of Steklov Mathematical Institute\newline
Fontanka 27, St.\,Petersburg \,191023, Russia}

\email{vt@pdmi.ras.ru}

\begin{document}

\maketitle

\section{Introduction}

The \KZv/ (\KZ/\)) \eq/s is a holonomic system of \difl/ \eq/s for correlation
\fn/s in conformal field theory on the sphere \cite{KZ}\). The \KZe/s play
an important role in \rep/ theory of affine Lie algebras and \qg/s, see for
example \cite{EFK}\). There are \rat/, \tri/ and elliptic versions of \KZe/s,
depending on what kind of coefficient \fn/s the \eq/s have. In this paper we
will consider only the \rat/ and \tri/ versions of the \KZe/s.
\vsk.1>
The \rat/ \KZe/s associated with a reductive Lie algebra $\gg$ is a system of
\eq/s for a \fn/ $u(\zn)$ of complex \var/s $\zn$, which takes values in
a tensor product $\Vox$ of \g-modules $\Vn$. The \eq/s depend on a complex
parameter $\ka$, and their coefficients are expressed in terms of the \sym/
tensor ${\Om\in U(\gg)\ox U(\gg)}$ corresponding to a nondegenerate invariant
bilinear form on $\gg$. For example, if $\gg=\slt$ and $e,f,\)h$ are
its standard generators \st/ $[\)e\>,\]f\>]=h$, then
${\Om=e\ox\}f+f\]\ox e+h\ox h/2}$.
\vsk.1>
The \rat/ \KZe/s are
\begin{gather}
\ka\>\frac{\der u}{\der z_i\]}\;=\,
\sum_{\tsize\fratop{j=1}{j\ne i}}^n\>\frac{\Om^{(ij)}}{z_i-z_j}\ u\,,
\qquad i=1\lc n\,,\kern-2em
\Tag{kz}
\cnn-.2>
\end{gather}
where ${\Om^{(ij)}\]\in\End(\Vox)}$ is the operator acting as $\Om$ on
${V_i\]\ox V_j}$ and as the identity on all other tensor factors; for instance,
\vvn.1>
$$
\Om^{(12)}(\vox)\)=\)\bigl(\Om\>(v_1\ox v_2)\bigr)\ox v_3\lox v_n\,.
$$
\par
All over the paper we will assume that $\ka$ is not a \rat/ number. Properties
of \sol/s of the \KZe/s depend much on whether $\ka$ is \rat/ or not.
\vsk.1>
Equations \(kz) can be generalized to a holonomic system of \difl/ \eq/s
depending on an element $\la\in\gg$:
\vvn.1>
\begin{gather}
\ka\>\frac{\der u}{\der z_i\]}\;=\,\la^{(i)}u\,+\)
\sum_{\tsize\fratop{j=1}{j\ne i}}^n\>\frac{\Om^{(ij)}}{z_i-z_j}\ u\,,
\qquad i=1\lc n\,.\kern-2em
\Tag{kza}
\cnn-.1>
\end{gather}
Here ${\la^{(i)}\]\in\End(\Vox)}$ acts as $\la$ on $V_i$ and as the identity
on all other tensor factors: \>$\la^{(i)}(\vox)\)=\)v_1\lox\la\)v_i\lox v_n$.
System \(kza) is also called the {\it\rat/ \>\KZe/s\/}.
\vsk.1>
Further on we will assume that $\la$ is a semisimple regular element of $\gg$.
Let $\hg\sub\gg$ be the Cartan subalgebra containing $\la$, and let
${e_\al\in\gg}$ be a root vector corresponding to a root $\al\in\hg^*\}$.
We normalize the root vectors by ${(e_\al\>,e_{-\al})=1}$, where $(\>{,}\>)$
is the bilinear form on $\gg$ corresponding to the tensor $\Om$.
\vsk.1>
In \cite{FMTV} system \(kza) was extended to a larger system of holonomic
\difl/ \eq/s for a \fn/ $u(\zn\);\la)$ on ${\Cn\}\oplus\hg}$. In addition to
\eq/s \(kza) the extended system includes the following \eq/s \wrt/ $\la$:
\vvn.1>
\begin{gather}
\ka\>\Dmu u\,=\,\sum_{i=1}^n\)z_i\)\mu^{(i)}u\,+\)\sum_\al
\frac{(\mu\),\al)}{2\>(\)\la\>,\al)}\;e_\al\)e_{-\al}\,u\,,
\qquad \mu\in\hg\,,\kern-2em
\Tag{dd}
\end{gather}
where $\Dmu$ is the directional derivative: $\Dmu
u(\)\la)\)=\)\bigl(\)\der_t\)u(\)\la+t\)\mu)\bigr)\big|_{t=0}$. Equations
\(dd) are called the {\it\rat/ dynamical \difl/\/} \)(\DD/\)) {\it\eq/s\/}.
\vsk.1>
A special case of \eq/s \(dd)\), when ${n=1}$ and ${z_1=0}$, was discovered
for a completely different reason. Around 1995 studying hyperplanes
arrangements De~Concini and Procesi introduced in an unpublished work
a connection on the set of regular elements of the Cartan subalgebra $\hg$.
The \eq/s for horizontal sections of the De~Concini\>-\)Procesi connection
coincide with the \rat/ \DD/ \eq/s. The same connection also appeared later
in \cite{TL}\). De~Concini and Procesi conjectured that the monodromy of their
connection is described in terms of the quantum Weyl group of type $\gg$.
For $\gg=\sln$ this conjecture was proved in \cite{TL}\).
\vsk.1>
If all \g-modules $\Vn$ are \hwm/s, \sol/s of the \KZe/s \(kz) can be written
in terms of multidimensional \hint/s \cite{SV}\), \cite{V}\). The construction
of \hgeom/ \sol/s can be generalized in a straightforward way to the case of
\KZe/s \(kza)\), see~\cite{FMTV}\). Moreover, it is shown in \cite{FMTV}
that the \hgeom/ \sol/s of the \KZe/s obey the \DD/ \eq/s \(dd) as well.
Generically, \hgeom/ \sol/s of the \KZ/ and \DD/ \eq/s are complete, that is,
they form a basis of \sol/s of those systems of \difl/ \eq/s.
\vsk.1>
An amusing fact about the \hgeom/ \sol/s is that though systems \(kza) and
\(dd) have rather similar look, and the \var/s $\zn$ and $\la$ seem to play
nearly interchangable roles, the formulae for the \hgeom/ \sol/s of the \KZ/
and \DD/ \eq/s involve $\zn$ and $\la$ in a highly non\sym/ way. While the
\var/s $\zn$ determine singularities of integrands of the \hint/s and enter
there in a rather complicated manner, $\la$ appears in the integrands only in
a very simple way via the exponential of a linear form. Such asymmetry suggests
the following idea. Suppose that a certain holonomic system of \difl/ \eq/s
can be viewed both as a special case of system \(kza) and as a special case of
system \(dd)\), maybe not for the same Lie algebra $\gg$. Then one can get two
types of integral formulae for \sol/s of that system, and \sol/s of one kind
should be linear combinations of \sol/s of the other kind. Thus, this can lead
to nontrivial relations between \hint/s of different dimensions.
\vsk.1>
It turns out that the mentioned idea indeed can be realized in the framework
of the $\glkn$ duality. This duality plays an important role in the \rep/
theory and the classical invariant theory, see \cite{Zh1}\), \cite{Ho}\).
It was observed in \cite{TL} that under the $\glkn$ duality the \KZe/s \(kz)
for the Lie algebra $\slk$ correspond to the \DD/ \eq/s \(dd) (with $n$
replaced by $k$ and all $z\}$'s being equal to zero) for the Lie algebra
$\sln$. This fact was used in \cite{TL} to compute the monodromy of the
De~Concini\>--\)Procesi connection in terms of the quantum Weyl group action.
\vsk.1>
Systems \(kza) and \(dd) are counterparts of each other under the $\glkn$
duality in general as well, see ~\cite{TV4}\). Employing this claim for
$k=n=2$, after all one arrives to identities for \hint/s of different
dimensions \cite{TV6}\). One can expect that there are similar identities
for \hint/s for an arbitrary pair $k\),n$.
\vsk.1>
There are various generalizations of the \KZe/s. The \fn/ $\Om/z$, describing
the coefficients of the \KZe/s, is the simplest example of a classical
\r-matrix --- a \sol/ of the classical \YB/. Starting from any classical
\r-matrix with a spectral parameter one can write down a holonomic system
of \difl/ \eq/s, see~\cite{Ch2}\). The obtained system is called the \KZe/s
associated with the given \r-matrix. For example, the standard \tri/
\r-matrix is
$$
r(z)\,=\;\frac\Om{z-1}\,+\,\frac12\>\sum_a\>\xi_a\]\ox\)\xi_a\,+\>
\sum_{\al>0}\>e_\al\]\ox\)e_{-\al}\,,
$$
where $\lb\>\xi_a\rb$ is an orthonormal basis of the Cartan subalgebra,
and the second sum is taken over all positive roots $\al$, \cf. \(rz) for
the Lie algebra $\glk$. The \tri/ \r-matrix satisfies the classical \YB/
\vvn.3>
$$
\bigl[\)r_{12}(z/w)\),r_{13}(z)+r_{23}(w)\bigr]\)+\)
\bigl[\)r_{13}(z)\),r_{23}(w)\bigr]\,=\,0\,.
\vv.3>
$$
The corresponding \KZe/s are
\begin{gather}
\ka\>z_i\>\frac{\der u}{\der z_i\]}\;=\,\la^{(i)}u\,+\)
\sum_{\tsize\fratop{j=1}{j\ne i}}^n\>r^{(ij)}(z_i/\)z_j)\;u\,,
\qquad i=1\lc n\,,\kern-2em
\Tag{kzt}
\end{gather}
where $\la$ is an element of the Cartan subalgebra. They are called the
{\it \tri/\/} \KZe/s associated with the Lie algebra $\gg$. System \(kza) can
\vvgood
be considered as a limiting case of system \(kzt) by the following procedure:
one replaces the \var/s $\zn$ by $e^{\)\eps z_1}\}\lc\)e^{\)\eps z_n}\}$ and
$\la$ by $\la/\eps$, and then sends $\eps$ to $0$.
\vsk.1>
The \dif/ analogue of the \KZe/s --- the {\it quantized \KZv/\/} (\qKZ/\))
{\it \eq/s\/} --- were introduced in \cite{FR}\). Coefficients of the \qKZe/s
are given in terms of quantum \R-matrices --- \sol/s of the quantum \YB/:
\vvn.3>
$$
R_{12}(z-w)\>R_{13}(z)\>R_{23}(w)\,=\,R_{23}(w)\>R_{13}(z)\>R_{12}(z-w)\,.
\vv.2>
$$
There are \rat/, \tri/ and elliptic versions of \KZe/s, the corresponding
\R-matrices coming from the \rep/ theory of Yangians, \qaff/ algebras and
elliptic \qg/s, \resp/. The \em{\rat/} \qKZe/s associated with the Lie algebra
$\gg$ is a holonomic system of \deq/s for a \fn/ $u(\zn)$ with values in
a tensor product ${\Vox}$ of modules over the Yangian $Y(\gg)$:
\vvn.2>
\begin{align}
u(z_1\lc z_i\]+\ka\lc z_n)\,=\,\bigl(& R_{1i}(z_1\]-z_i\]-\ka)\ldots
R_{i-1,i}(z_{i-1}\]-z_i\]-\ka)\bigr)\1\x{}\kern-1em
\Tag{qkz}
\\[4pt]
{}\x\,(\emu)^{\](i)}\>&
R_{in}(z_i\]-z_n)\ldots R_{i,i+1}(z_i\]-z_{i+1})\,u(\zn)\,,
\kern-1em
\notag
\cnn.2>
\end{align}
$i=1\lc n$. Here $\mu$ is an element of the Cartan subalgebra and $R_{ij}(z)$
is the \R-matrix for the tensor product $V_i\ox V_j$ of the Yangian modules.
\vsk.1>
There are also several generalizations of the \rat/ \difl/ dynamical \eq/s.
The \dif/ analogue of the \DD/ \eq/s --- the {\it \rat/ \dif/ dynamical\/}
(\qDD/\)) {\it\eq/s\/} --- was suggested in \cite{TV3}\). The idea was to
extend the \tri/ \KZe/s \(kzt) by \eq/s \wrt/ $\la$ similarly to the way in
which system \(dd) extends the \rat/ \KZe/s \(kza)\), and to obtain a holonomic
system of \difl/\>-\)\deq/s for a \fn/ $u(\zn\);\la)$ on $\Cn\}\oplus\hg$.
The \rat/ \qDD/ \eq/s have the form
\vvn.3>
\begin{gather}
u(\zn\);\la+\ka\)\om)\,=\,Y_\om(\zn\);\la)\>u(\zn\);\la)
\Tag{qdd}
\cnn.2>
\end{gather}
where $\om$ is an integral \wt/ of $\gg$, and the operators $Y_\om$ are
written in terms of the extremal cocycle on the Weyl group of $\gg$.
The extremal cocycles and their special values, the extremal projectors,
are important objects in the \rep/ theory of Lie algebras and Lie groups,
see~\cite{AST}\), \cite{Zh2}\), \cite{Zh3}\), \cite{ST}\).
\vsk.1>
The ideas used in \cite{TV3} were further developed in \cite{EV} where a new
concept of the dynamical Weyl group was introduced, and the \tri/ version of
the \dif/ dynamical \eq/s was suggested.
\vsk.1>
There is also the \tri/ version of the \difl/ dynamical \eq/s, which,
in principle, can be obtained by degenerating the \tri/ \dif/ dynamical
\eq/s. The explicit form of the \tri/ \difl/ dynamical \eq/s for the Lie
algebras $\glk$ and $\slk$ was obtained in \cite{TV4} by extending the \rat/
\qKZe/s \(qkz) by \eq/s \wrt/ $\mu$ in such a way that the result is
a holonomic system of \dif/\>-\)\difl/ \eq/s for a \fn/ $u(\zn\);\mu)$
on $\Cn\}\oplus\hg$.
\vsk.1>
The $\glkn$ duality naturally applies to the \tri/ and \dif/ versions of the
\KZ/ and dynamical \eq/s. Under the duality, the \tri/ \KZe/s \(kzt) for the
Lie algebra $\glk$ correspond to the \tri/ \difl/ dynamical \eq/s for the Lie
algebra $\gln$, and vice versa. At the same time the \rat/ \qKZe/s for $\glk$
are counterparts of the \rat/ \qDD/ \eq/s for $\gln$. To relate the \tri/ \qKZ/
and \qDD/ \eq/s, one has to employ the \$q$-analogue of the $\glkn$ duality:
the $\Uqkn$ duality described in \cite{B}\), \cite{TL}\).
\vsk.1>
Hypergeometric \sol/s of the \tri/ \KZe/s \(kzt) can be written almost
in the same manner as those of the \rat/ \KZe/s \(kz)\), see~\cite{Ch1}\),
\cite{MV}\). Conjecturally, the \hgeom/ \sol/s of the \tri/ \KZe/s obey
the corresponding \rat/ \qDD/ \eq/s. For the Lie algebra $\slk$ this claim
was proved in \cite{MV}\). On the other hand, \sol/s of the \rat/ \qKZe/s
can be written in terms of suitable \q-\hgeom/ Jackson integrals \cite{TV1}\),
or \q-\hint/s of Mellin\)-Barnes type \cite{TV2}\). Thus, using the $\glkn$
duality, one can obtain \sol/s of a certain system of \difl/\>-\)\deq/s
both in terms of ordinary \hint/s and \q-\hint/s of Mellin\)-Barnes type,
and establish nontrivial relations between those integrals. For $k=n=2$ this
has been done in \cite{TV7}\). The obtained relations are multidimensional
analogues of the equality of two integral \rep/s for the Gauss \hgeom/ \fn/
$_2F_1$:
\vvn.3>
\begin{align*}
{}_2F_1(\al\),\bt\);\gm\);z)\,
& {}=\;\frac{\Gm(\gm)}{\Gm(\al)\>\Gm(\gm-\al)}
\,\,\int_{\!0}^{1\!}u^{\al-1}\>(1-u)^{\gm-\al-1}\>(1-u\)z)^{-\)\bt}\>du
\\[3pt]
& {}=\;\frac1{2\pi i}\;\frac{\Gm(\gm)}{\Gm(\al)\>\Gm(\bt)}
\,\int_{-\)i\)\infty\)-\eps}^{+\)i\)\infty\)-\)\eps}\!(-\)z)^s\;
\frac{\Gm(-\)s)\>\Gm(s+\al)\>\Gm(s+\bt)}{\Gm(s+\gm)}\ ds\,.
\end{align*}
\vsk.3>
As it was pointed out by J\&Harnad, the duality between the \KZ/ and \DD/ \eq/s
in the \rat/ \difl/ case is essentially the ``quantum'' version of the duality
for isomonodromic deformation systems \cite{H1}\). The relation of the \difl/
\KZe/s and the isomonodromic deformation systems is described in \cite{R}\),
\cite{H2}\). From this point of view the \rat/ \qDD/ \eq/s can be considered
as ``quantum'' analogues of the Schlesinger transformations, though
the correspondence is not quite straightforward.
\vsk.2>
The paper is organized as follows. After introducing basic notation we
subsequently describe the \difl/ \KZ/ and \DD/ \eq/s, and the \rat/ \dif/ \qKZ/
and \qDD/ \eq/s, for the Lie algebra $\glk$. This is done in Sections~2\,--\,5.
Then we consider the $\glkn$ duality in application to the \KZ/ and dynamical
\eq/s. In the last two sections we describe the \hgeom/ \sol/s of the \eq/s,
and use the duality relations to establish identities for \hgeom/ and
\q-\hint/s of different dimensions.

\section{Basic notation}
Let $n$ be a nonnegative integer. A partition $\la=(\la_1\),\la_2\),\ldots{})$
with at most $\>k$ parts is an infinite nonincreasing sequence of nonnegative
integers \st/ $\la_{k+1}=0$. Denote by $\Pck$ the set of partitions with at
most $k$ parts and by $\Pc$ the set of all partitions. We often make use of the
embedding $\Pck\]\to\Ck$ given by truncating the zero tail of a partition:
$(\la_1\)\lc\la_k\),0\>,0\>,\ldots{})\map(\la_1\)\lc\la_k)$. Since obviously
$\Pcm\]\sub\Pck$ for $m\le k$, in fact, one has a collection of embeddings
$\Pcm\to\Ck$ for any $m\le k$. What particular embedding is used will be clear
from the context.
\vsk.1>
Let $e_{ab}$, $a,b=1\lc k$, be the standard basis of the Lie algebra $\glk$:
$[\)e_{ab}\>,e_{cd}\)]\)=\)\dl_{bc}\>e_{ad}-\dl_{ad}\>e_{cb}$. We take
the Cartan subalgebra $\hg\sub\glk$ spanned by $e_{11}\lc e_{kk}$, and
the nilpotent subalgebras $\ngp$ and $\ngm$ spanned by the elements $e_{ab}$
for ${a<b}$ and ${a>b}$, \resp/. One has the standard Gauss decomposition
$\glk=\ngp\]\oplus\)\hg\oplus\ngm$.
\vsk.1>
Let $\eps_1\lc\eps_k$ be the basis of $\)\hg^*\}$ dual to $e_{11}\lc e_{kk}$:
$\bra\)\eps_a\>,e_{bb}\)\ket\)=\)\dl_{ab}$. We identify $\hg^*\}$ with $\Ck\}$
mapping $\la_1\)\eps_1\lsym+\la_k\>\eps_k$ to $(\la_1\lc\la_k)$.
The root vectors of $\glk$ are $e_{ab}$ for $a\ne b$, the corresponding root
being equal to $\al_{ab}=\)\eps_a\]-\eps_b$. The roots $\al_{ab}$ for $a<b$
are positive.
\vsk.1>
We choose the standard invariant bilinear form $(\,{,}\,)$ on $\glk$:
$(\)e_{ab}\>,e_{cd}\))\)=\)\dl_{ad}\>\dl_{bc}$. It defines an \iso/
$\hg\to\hg^*\}$. The induced bilinear form on $\hg^*$ is
$(\)\eps_a\>,\eps_b\))\)=\)\dl_{ab}$.
\vsk.1>
\vsk.1>
For a \$\glk$-module $W\}$ and a \wt/ $\la\in\hg^*\}$
let $\Wla$ be the \wt/ subspace of $W\}$ of \wt/ $\la$.
\vsk.2>
For any $\la\in\Pck$ we denote by $V_\la$ the \irr/ \$\glk$-module with
\hw/ $\la$. By abuse of notation, for any $l\in\Zp$ we write $V_l$ instead
of $V_{(\)l,\)0\)\lc 0)}$. Thus, $V_0=\C$ is the trivial \$\glk$-module,
$V_1=\Ck$ with the natural action of $\glk$, and $V_l$ is the \$l$-th symmetric
power of $V_1$.
\vsk.1>
Define a \$\glk$-action on the \pol/ ring $\C\)[\)x_1\lc x_k]$ by \difl/
operators: $e_{ab}\map x_a\der_b$, where $\der_b=\der/\der x_b$, and denote
the obtained \$\glk$-module by $\VV\]$. Then
\vvn-.3>
\begin{gather}
\VV\>=\,\Plus_{l=0}^\infty\,V_l\,,
\Tag{VV}
\cnn.2>
\end{gather}
the submodule $V_l$ being spanned by homogeneous \pol/s of degree $l$.
The \hwv/ of the submodule $V_l$ is $x_1^{\>l}$.

\section{\KZv/ and \difl/ dynamical \eq/s}
For any $g\in\Uk$ set
$g^{(i)}=\>\id\)\lox{}\%{g}_{\Clap{\sss\text{\=$i$-th}}}{}\]\lox\)\id
\in\bigl(\Uk\bigr)\vpb{\ox n}\!$.
\vv-.2>
We consider $\Uk$ as a subalgebra of $\Ukn\}$, the embedding \>$\Uk\)\hto\Ukn$
\vvn.06>
being given by the \$n\)$-fold coproduct, that is,
$x\>\map x^{(1)}\lsym+x^{(n)}$ for any $x\in\glk$.
\vsk.1>
Let \,$\Om\,=\}\sum_{a,b=1}^k e_{ab}\)\ox\)e_{ba}$ be the Casimir tensor,
and let
\vvn-.2>
\begin{gather*}
\Om_+\>=\;\frac12\>\sum_{a=1}^k e_{aa}\)\ox\)e_{aa}\>+\!
\sum_{1\le a<b\le k}\!\}e_{ab}\)\ox\)e_{ba}\,,
\\[4pt]
\Om_-\>=\;\frac12\>\sum_{a=1}^k e_{aa}\)\ox\)e_{aa}\>+\!
\sum_{1\le a<b\le k}\!\}e_{ba}\)\ox\)e_{ab}\,,
\cnn.2>
\end{gather*}
so that ${\Om=\Om_+\}+\)\Om_-}$.
The standard \tri/ \r-matrix, associated with the Lie algebra $\glk$ is
\vvn.2>
\begin{gather}
r(z)\,=\;\frac\Om{z-1}\,+\,\Om_+\,=\;\frac{z\>\Om_+\]+\>\Om_-}{z-1}\;.
\Tag{rz}
\cnn.2>
\end{gather}
\par
Fix a nonzero complex number $\ka$. Consider \difl/ operators $\nablasn\]$ and
\>$\hablasn\]$ with coefficients in $\Ukn\}$ depending on complex \var/s $\zn$
and $\lak$:
\vvn.2>
\begin{gather}
\nabla_{\}z_i}(z\);\la)\,=\,\ka\>\frac\der{\der z_i\]}
\>-\>\sum_{a=1}^k\>\la_a\>({e_{aa}})^{(i)}
\>-\>\sum_{\tsize\fratop{j=1}{j\ne i}}^n\>\frac{\Om^{(ij)}}{z_i-z_j}\;,
\Tag{nabla}
\\[4pt]
\habla_{\}z_i}(z\);\la)\,=\,\ka\)z_i\>\frac\der{\der z_i\]}\>-\>
\sum_{a=1}^k\>\bigl(\la_a\]-\frac{e_{aa}}2\bigr)\>({e_{aa}})^{(i)}
\>-\>\sum_{\tsize\fratop{j=1}{j\ne i}}^n\>r^{(ij)}(z_i/\)z_j)\,.
\Tag{habla}
\cnn.1>
\end{gather}
The \difl/ operators $\nablasn$ (resp.~$\hablasn$\}) are called
the \em{\rat/} (resp.~\em{\tri/}) \em{\KZv/} (\KZ/\)) \em{operators}.
The following statements are well known.
\begin{theorem}
\label{ratKZ}
The operators $\>\nablasn\}$ pairwise commute.
\end{theorem}
\begin{theorem}
\label{trigKZ}
The operators $\>\hablasn\}$ pairwise commute.
\end{theorem}
The \rat/ \KZe/s associated with the Lie algebra $\glk$ is a system
of \difl/ \eq/s
\begin{gather}
\nabla_{\}z_i}u\,=\,0\,,\qqq i=1\lc n\,,\kern-2em
\Tag{KZ}
\cnn.1>
\end{gather}
for a \fn/ $u(\zn\);\lak)$ taking values in an \$n\)$-fold tensor product of
\$\glk$-modules. Similarly, the \tri/ \KZe/s associated with the Lie
algebra $\glk$ is a system of \difl/ \eq/s
\vvn.5>
\begin{gather}
\habla_{\}z_i}u\,=\,0\,,\qqq i=1\lc n\,.\kern-2em
\Tag{KZt}
\cnn.3>
\end{gather}
for a \fn/ $u(\zn\);\lak)$.
\par
Introduce \difl/ operators $\Dlk\]$ and \>$\Dhk\]$ with coefficients
in $\Ukn\}$ depending on complex \var/s $\zn$ and $\lak$:
\vvn-.3>
\begin{align}
D_{\la_a}(z\);\la)\, &{}=\,\ka\>\frac\der{\der\la_a\]}
\>-\)\sum_{i=1}^n\)z_i\>({e_{aa}})^{(i)}
\>-\)\sum_{\tsize\fratop{b=1}{b\ne a}}^k\>
\frac{e_{ab}\>e_{ba}-\)e_{aa}}{\la_a\]-\la_b}\;.
\Tag{D}
\\[3pt]
\Dh_{\la_a}(z\);\la)\, &{}=\,\ka\>\la_a\frac\der{\der \la_a\]}\>+
\>\frac{e_{aa}^{\)2}}2\>-\)\sum_{i=1}^n\)z_i\>({e_{aa}})^{(i)}\>-{}
\\[3pt]
&{}-\)\sum_{b=1}^k\,\sum_{1\le i<j\le n}\!(e_{ab})^{(i)}\)(e_{ba})^{(j)}
\>-\)\sum_{\tsize\fratop{b=1}{b\ne a}}^k\)
\frac{\la_b}{\la_a\]-\la_b\]}\,(e_{ab}\>e_{ba}-\)e_{aa})\,.\kern-1em
\notag
\cnn-.2>
\end{align}
Recall that $e_{ab}\)=\sum_{i=1}^n\>(e_{ab})^{(i)}$. The operators $\Dlk$
(resp.~$\Dhk$) are called the \em{\rat/} (resp.~\em{\tri/}) \em{\difl/
dynamical} (\)\DD/\)) \em{operators}.
\begin{theorem}
\label{KZDD}
The operators $\>\nablasn$, $\Dlk\}$ pairwise commute.
\end{theorem}
\nt
The theorem follows from the same result for the \rat/ \KZ/ and \DD/ operators
associated with the Lie algebra $\slk$, see~\cite{FMTV}\).
\begin{theorem}
\label{trigDD}
\cite{TV4}
The operators $\>\Dhk\}$ pairwise commute.
\end{theorem}
\nt
The statement can be verified in a straightforward way.
\vsk.5>
Later we will formulate analogues of Theorem~\ref{KZDD} for the \tri/ \KZ/
operators and the \tri/ \DD/ operators, see Theorems~\ref{KZqDD} and
\ref{qKZDD}\). They involve \dif/ dynamical operators and \dif/ (quantized)
\KZv/ operators which are discussed in the next two sections.
\vsk.2>
The \rat/ \DD/ \eq/s associated with the Lie algebra $\glk$ is a system
of \difl/ \eq/s
\vvn.1>
\begin{gather}
D_{\la_a}u\,=\,0\,,\qqq a=1\lc k\,,\kern-2em
\Tag{DD}
\cnn.1>
\end{gather}
for a \fn/ $u(\zn\);\lak)$ taking values in an \$n\)$-fold tensor product of
\$\glk$-modules. Similarly, the \tri/ \DD/ \eq/s associated with the Lie
algebra $\glk$ is a system of \difl/ \eq/s
\vvn.5>
\begin{gather}
\Dh_{\la_a}u\,=\,0\,,\qqq a=1\lc k\,.\kern-2em
\Tag{DDt}
\cnn.2>
\end{gather}
for a \fn/ $u(\zn\);\lak)$.
\begin{remark}
Systems \(KZt) and \(DD) are not precisely the same as specializations of
the respective systems \(kzt) and \(dd) for the Lie algebra $\glk$. However,
in both cases the \dif/ is not quite essential and can be worked out. The form
of the operators $\habla_{\}z_i\}}$ and $D_{\la_a}$ given in this section,
see~\(habla) and \(D), fits the best the framework of the $\glkn$ duality.
\end{remark}

\section{Rational \dif/ dynamical \eq/s}
For any $a\>,b=1\lc k$, $a\ne b$, \>introduce a series $B_{ab}(t)$ depending
on a complex \var/ $t$:
\vvn.1>
$$
B_{ab}(t)\,=\,1+\>\sum_{s=1}^\infty\>e_{ba}^s\)e_{ab}^s\,
\prod_{j=1}^s\,\frac1{j\>(t-\)e_{aa}\]+\)e_{bb}-j)}\;.
\vv.3>
$$
The series has a well-defined action on any \fd/ \$\glk$-module $W\}$, giving
an \$\End(W)\)$-valued \raf/ of $t$. The series $B_{ab}(t)$ have zero \wt/:
\begin{gather}
\bigl[\)B_{ab}(t)\>,x\)\bigr]\>=\>0\qquad\text{for any}\quad x\in\hg\,,
\Tag{Bx}
\cnn.3>
\end{gather}
satisfy the inversion relation
\vvn.2>
\begin{gather}
B_{ab}(t)\>B_{ba}(-\)t)\,=\,1\,-\frac{e_{aa}\]-e_{bb}}t\;,
\Tag{BB}
\cnn.1>
\end{gather}
and the braid relation
\vvn.4>
\begin{gather}
B_{ab}(t-s)\>B_{ac}(t)\>B_{bc}(s)\,=\,B_{bc}(s)\>B_{ac}(t)\>B_{ab}(t-s)\,.
\Tag{BBB}
\cnn.3>
\end{gather}
Relation \(Bx) is clear. Relations \(BB) and \(BBB) follow from \cite{TV3}\),
namely from the properties of \fn/s $B_w(\la)$ considered there in the $\slk$
case, see~\cite[Section~2.6\)]{TV1}\). In notation of \cite{TV3} the series
$B_{ab}(t)$ equals ${p\)(t-1\);\)e_{aa}\]-e_{bb}\),e_{ab}\),e_{ba})}$.
\begin{remark}
The series $B_{ab}(t)$ first appeared in the definition of the extremal
projectors \cite{AST} and the extremal cocycles on the Weyl group
\cite{Zh2}\), \cite{Zh3}\).
\end{remark}
Consider the products $X_1\lc X_k$ depending on complex \var/s $\zn$ and
$\lak$:
\vvn.2>
\begin{align}
X_a(z\);\la)\,= {}& \,\bigl(\)B_{ak}(\la_{ak})\ldots
B_{a,\)a+1}(\la_{a,\)a+1})\)\bigr)\vpb{-1}\x{}
\Tag{Xa}
\\[4pt]
& {}\x\,\prod_{i=1}^n\,\bigl(z_i^{-\)e_{aa}}\bigr)\'i\,
B_{1a}(\la_{1a}\]-\ka)\ldots B_{a-1,\)a}(\la_{a-1,\)a}\]-\ka)\,,\kern-2em
\notag
\cnn.3>
\end{align}
where $\la_{bc}=\la_b-\la_c$. They act on any \$n\)$-fold tensor product
${W_1\lox W_n}$ of \fd/ (more generally, \hw/) \$\glk$-modules.
\vsk.3>
Let $T_u$ be a \dif/ operator acting on a \fn/ $f(u)$ by the rule
\vvn.4>
$$
(T_uf)(u)\,=\,f(u+\ka)\,.
\vv.3>
$$
Introduce \dif/ operators $\Qlk$:
\vvn.4>
$$
Q_{\la_a}(z\);\la)\,=\,X_a(z\);\la)\,T_{\la_a}\,.
\vv.3>
$$
They are called the \em{\rat/ \dif/ dynamical\)} (\qDD/\)) operators.
\begin{theorem}
\label{KZqDD}
The operators $\>\hablasn$, $\Qlk\}$ pairwise commute.
\end{theorem}
\nt
The theorem follows from the same result for the \tri/ \KZ/ and \rat/ \qDD/
operators in the $\slk$ case, see~\cite{TV1}\). Theorem~\ref{KZqDD} extends
Theorem~\ref{trigKZ}\), and is analogous to Theorem~\ref{KZDD}\).
\vsk.5>
In more conventional form the equalities
\vvn.3>
$$
[\)\habla_{\}z_i}\),\)Q_{\la_a}\)]\)=\)0\,,\qqq
[\)Q_{\la_a}\),\)Q_{\la_b}\)]\)=\)0\,,
\vv.3>
$$
\resp/ look like
\vvn.4>
\begin{gather*}
\habla_{\}z_i}(z\);\la)\,X_a(z\);\la)\,=\,
X_a(z\);\la)\,\habla_{\}z_i}(z\);\la_1\lc\la_a\]+\ka\lc \la_k)\,,
\\[7pt]
X_a(z\);\la)\,X_b(z\);\la_1\lc\la_a\]+\ka\lc \la_k)\,=\,
X_b(z\);\la)\,X_a(z\);\la_1\lc\la_b\]+\ka\lc \la_k)\,.
\end{gather*}
\vsk.5>
The \em{\rat/ \dif/ dynamical} (\qDD/\)) \em{\eq/s\,} associated with the Lie
algebra $\glk$ is a system of \deq/s
\vvn.4>
\begin{gather}
Q_{\la_a}u\,=\,u\,,\qqq a=1\lc k\,,\kern-2em
\Tag{qDD}
\cnn.2>
\end{gather}
for a \fn/ $u(\zn\);\lak)$ taking values in an \$n\)$-fold tensor product
of \$\glk$-modules.

\section{Rational \dif/ \KZv/ \eq/s}

For any two \irr/ \fd/ \$\glk$-modules $V,\)W\}$ there exists a distinguished
\$\End(V\]\ox W)\)$-valued \raf/ $\RVW(t)$ called the \em{rational \R-mat\-rix}
for the tensor product $V\]\ox W\}$. The definition of $\RVW(t)$ comes from
the \rep/ theory of the Yangian $\Yk$.
\vsk.1>
The Yangian $\Yk$ is an in\fd/ Hopf algebra, which is a flat deformation
of the universal enveloping algebra $\Ukx$ of \$\glk$-valued \pol/ \fn/s.
The subalgebra of constant \fn/s in $\Ukx$, which is isomorphic to $\Uk$, is
preserved under the deformation. Thus, the algebra $\Uk$ is embedded in $\Yk$
as a Hopf subalgebra, and we identify $\Uk$ with the image of this embedding.
\vsk.1>
There is an algebra \hom/ $\ev:\)\Yk\)\to\>\Uk$, \>called the \em{evaluation
\hom/}, which is identical on the subalgebra $\Uk\sub\Yk$. It is a deformation
of the \hom/ $\Ukx\)\to\>\Uk$ which sends any \pol/ to its value at $x=0$.
The evaluation \hom/ is not a \hom/ of Hopf algebras.
\vsk.2>
The Yangian $\Yk$ has a distinguished one-parametric family of \aut/s
$\rho_u$ depending on a complex parameter $u$, which is informally called
the \em{shift of the spectral parameter}. The \aut/ $\rho_u$ corresponds to
the \aut/ $p(x)\)\map p(x+u)$ of the Lie algebra $\glk[\)x\)]$.
For any \$\glk$-module $W\}$ we denote by $W(u)$ the pullback of $W\}$
via the \hom/ \>$\ev\o\,\rho_u$. Yangian modules of this form are called
\em{evaluation modules}.
\vsk.2>
For any \fd/ \irr/ \$\glk$-modules $V\),\>W\}$ the tensor products
$V(t)\ox W(u)$ and $W(u)\ox V(t)$ are isomorphic \irr/ \$\Yk$-modules,
provided $t-u\nin\Z$. The intertwiner $V(t)\ox W(u)\)\to\)W(u)\ox V(t)$
can be taken of the form $\PVW\>\RVW(t-u)$, where
\>$\PVW\}:\)V\]\ox W\]\to\)W\]\ox V$ is the flip map:
$\PVW\}:\)v\ox w\)\map\alb\)w\ox v$, \,and $\RVW(t)$ is a \rat/
\$\End(V\]\ox W)\)$-valued \fn/, the \em{\rat/ \Rm/} for the tensor product
$V\ox W\}$.
\vsk.1>
The \Rm/ $\RVW(t)$ can be described in terms of the $\glk$ actions on
the spaces $V\}$ and $W\}$. It is determined uniquely up to a scalar multiple
by the \$\glk$-invariance,
\vvn.5>
\begin{gather}
\bigr[\)\RVW(t)\>,\>g\otimes 1+1\otimes g\>\bigr]\,=\,0\qqq
\text{for any}\quad g\in\glk\,,\kern-1.4em
\Tag{Rinv}
\cnn.4>
\end{gather}
and the commutation relations
\vvn.3>
\begin{gather}
\RVW(t)\>\Bigl(\)t\>e_{ab}\ox 1\,+\sum_{c=1}^k\>e_{ac}\ox e_{cb}\)\Bigr)\,=\,
\Bigl(\)t\>e_{ab}\ox 1\,+\sum_{c=1}^k\>e_{cb}\ox e_{ac}\)\Bigr)\>\RVW(t)\,.
\kern-1em
\Tag{Rdef}
\cnn.2>
\end{gather}
The standard normalization condition for $\RVW(t)$ is to preserve
the tensor product of the respective \hwv/s $v,w$:
\vvn.4>
$$
\RVW(t)\,v\ox w\,=\,v\ox w\,.
\vv.4>
$$
The introduced \Rms/ obey the inversion relation
\vvn.4>
\begin{gather}
\RVW(t)\,R_{WV}^{(21)}(-\)t)\,=\,1\,,
\Tag{inv}
\cnn.2>
\end{gather}
where $R_{WV}^{(21)}\)=\)P_{WV}\:\>R_{WV}\:\>\PVW$, \,and the \YB/
\vvn.4>
\begin{gather}
\RUV(t-u)\>\RUW(t)\>\RVW(u)\,=\,\RVW(u)\>\RUW(t)\>\RUV(t-u)\,.\kern-1em
\Tag{YB}
\end{gather}
\vsk.5>
The aforementioned facts on the Yangian $\Yk$ are well known.
A good introduction into the \rep/ theory of the Yangian $\Yk$ can be found
in \cite{MNO}\).
\vsk.2>
Consider the \$\glk$-module $\VV\]$, and let $V_l\sub\VV\]$ be the \irr/
component with \hwv/ $x_1^{\>l}$, see~\(VV)\). We define the \Rm/ $\RVV(t)$
to be a direct sum of the \Rms/ $\RVlm(t)$:
\vvn.4>
$$
\RVV(t)\,v\ox v'\)=\,\RVlm(t)\,v\ox v'\qqq
v\in\)V_l,\ \;v'\]\in\)V_m\,.\kern-2em
\vv.3>
$$
It is clear that $\RVV(t)$ obeys relations \(Rinv) and \(Rdef),
as well as the inversion relation and the \YB/.
\vsk.2>
Consider the products $K_1\lc K_n$ depending on complex \var/s $\zn$ and
$\lak$:
\vvn.2>
\begin{align}
\kern.8em K_i(z\);\la)\,= {}& \,\bigl(\)R_{in}(z_{in})\ldots
R_{i,\)i+1}(z_{i,\)i+1})\)\bigr)\vpb{-1}\x{}
\Tag{Ki}
\\[4pt]
&{}\x\,\prod_{a=1}^k\,\bigl(\)\Rlap{\la_a}{\la}\vpb{-\)e_{aa}}\bigr)\vpb{(i)}\>
R_{1i}(z_{1i}\]-\ka)\ldots R_{i-1,\)i}(z_{i-1,\)i}\]-\ka)\,,\kern-1.4em
\notag
\cnn.3>
\end{align}
acting on a tensor product ${W_1\lox W_n}$ of \$\glk$-modules.
Here $z_{ij}=z_i-z_j$, \,and
${R_{ij}(t)\)=\)\bigl(R_{W_iW_j}\:(t)\]\bigr)\vpb{(ij)}}$.
\vsk.2>
Introduce \dif/ operators $\Zzn$:
\vvn.4>
$$
Z_{z_i}(z\);\la)\,=\,K_i(z\);\la)\,T_{z_i}\,.
\vv.4>
$$
They are called the \em{\rat/ quantized \KZv/} (\qKZ/\)) operators. The next
theorem extends Theorem \ref{trigDD} and is analogous to Theorem~\ref{KZDD}\).
\begin{theorem}
\cite{FR}\), \cite{TV4}
\label{qKZDD}
The operators $\>\Zzn$, $\Dhk\}$ pairwise commute.
\end{theorem}
\nt
The \qKZ/ operators $\Zzn\}$ were introduced in \cite{FR}\), and their
commutativity was established therein. The fact that the \qKZ/ operators
commute with the operators $\Dhk\}$ can be verified in a straightforward way
using relations \(Rinv) and \(Rdef) for the \R-matrices.
\vsk.5>
In more conventional form the equalities
\,$[\)Z_{z_i}\),\)Z_{z_j}\)]\)=\)0$ \,and
\,$[\)Z_{z_i}\),\)\Dh_{\la_a}\)]\)=\)0$ \resp/ look like:
\vvn.5>
\begin{gather*}
K_i(z\);\la)\,K_j(z_1\lc z_i+\ka\lc z_n\);\la)\,=\,
K_j(z\);\la)\,K_i(z_1\lc z_j+\ka\lc z_n\);\la)\,,
\\[7pt]
\Dh_{\la_a}(z\);\la)\,K_i(z\);\la)\,=\,
K_i(z\);\la)\,\Dh_{\la_a}(z_1\lc z_i+\ka\lc z_n\);\la)\,.
\end{gather*}
\vsk.5>
The \em{\rat/ \qKZe/s\,} associated with the Lie algebra $\glk$
is a system of \deq/s
\vvn.2>
\begin{gather}
Z_{z_i}u\,=\,u\,,\qqq i=1\lc n\,,\kern-2em
\Tag{qKZ}
\cnn.2>
\end{gather}
for a \fn/ $u(\zn\);\lak)$ taking values in an \$n\)$-fold tensor product
of \$\glk$-modules.

\section{\boldmath$\{\glkn$ duality}
In this section we are going to consider the Lie algebras $\glk$ and $\gln$
simultaneously. In order to distinguish generators, modules, etc., we will
indicated the dependence on $k$ and $n$ explicitly, for example, $e_{ab}\+k\>$,
$V_\la\+n\}$.
\vsk.2>
Consider the \pol/ ring ${\PP_{kn}\)=\)\C\)[\)\xkn\)]}$ of $\>k\)n$ \var/s.
There are two natural \iso/s of vector spaces:
\vvn.4>
\begin{align}
\bigl(\C\)[\)x_1\lc x_k]\)\bigr)\vpb{\ox n}\} &{}\to\,\PP_{kn}\,,
\Tag{kiso}
\\[3pt]
(p_1\lox p_n)\)(\)x_{11}\lc x_{kn})\,
&{}=\,\prod_{i=1}^n\,p_i(\)x_{1i}\lc x_{ki})\,,
\notag
\\[-2pt]
\Text{and}
\bigl(\C\)[\)x_1\lc x_n]\)\bigr)\vpb{\ox k}\} &{}\to\,\PP_{kn}\,,
\Tag{niso}
\\[3pt]
(p_1\lox p_k)\)(\)x_{11}\lc x_{kn})\,
&{}=\,\prod_{a=1}^k\,p_a(\)x_{a1}\lc x_{an})\,.
\notag
\end{align}
\vsk.2>
Define a \$\glk$-action on $\PP_{kn}$ by
\vvn.1>
\begin{gather}
e_{ab}\+k\,\map\,\sum_{i=1}^n\,x_{ai}\)\der_{b\)i}\,,
\Tag{left}
\cnn.1>
\end{gather}
where $\der_{b\)i}=\der/\der x_{b\)i}$, and a \$\gln$-action by
\vvn.2>
\begin{gather}
e_{ij}\+n\,\map\,\sum_{a=1}^k\,x_{ai}\)\der_{aj}\,.
\Tag{right}
\end{gather}
\begin{proposition}
\label{Pkn}
As a \$\glk$-module, $\PP_{kn}$ is isomorphic to $\VVkn\!$ by \(kiso).
As a \$\gln$-module, $\PP_{kn}$ is isomorphic to $\VVnk\!$ by \(niso).
\end{proposition}
It is easy to see that the actions \(left) and \(right) commute with
each other, thus making $\PP_{kn}$ into a module over the direct sum
$\glk\]\oplus\)\gln$. The following theorem is well known.
\begin{theorem}
\label{dual}
The $\glk\]\oplus\)\gln$ module $\PP_{kn}$ has the decomposition
\vvn.4>
$$
\PP_{kn}\,=\,\Plus_{\la\)\in\)\Pc\_{\]\min(k,n)}}\ V_\la\+k\]\ox\)V_\la\+n\,.
$$
\end{theorem}
The module $\PP_{kn}$ plays an important role in the \rep/ theory and
the classical invariant theory, see \cite{Zh1}\), \cite{Ho}\), \cite{N}\).
\vsk.2>
Consider the action of \KZ/, \qKZ/, \DD/ and \qDD/ operators for the Lie
algebras $\glk$ and $\gln$ on \$\PP_{kn}$-valued \fn/s of $\zn$ and $\lak$,
treating the space $\PP_{kn}$ as a tensor product $\VVkn\!$ of \$\glk$-modules,
\vv.06>
and as a tensor product $\VVnk\!$ of \$\gln$-modules. If $F$ and $G$ act
on the \$\PP_{kn}$-valued \fn/s in the same way, we will write $F\simeq G$.
For instance,
\>${\bigl(e\+k_{aa}\bigr)\vpb{(i)}\]\simeq\)\bigl(e\+n_{ii}\bigr)\vpb{(a)}}$
since both $\bigl(e\+k_{aa}\bigr)\vpb{(i)}\}$ and
$\bigl(e\+k_{aa}\bigr)\vpb{(i)}\}$ act on $\PP_{kn}$ as $x_{ai}\)\der_{ai}$.
\vsk.2>
Introduce the following operators:
\vvn.4>
\begin{align}
C\+k_{ab}(t)\, &{}=\;\frac{\Gm(t+1)\,\Gm(t-e\+k_{aa}\]+e\+k_{bb})}
{\Gm(t-e\+k_{aa})\,\Gm(t+e\+k_{bb}+1)}\;,
\Tag{C}
\\[6pt]
C\+n_{ij}(t)\, &{}=\;\frac{\Gm(t+1)\,\Gm(t-e\+n_{ii}\]+e\+n_{jj})}
{\Gm(t-e\+n_{ii})\,\Gm(t+e\+n_{jj}+1)}\;.\kern-1em
\notag
\end{align}
\vsk.2>
\begin{theorem}
\cite{TV4}
\label{dualKZDD}
For any \,$i=1\lc n$ \>and \>$a=1\lc k$ \>we have
\vvn.6>
\begin{alignat}2
& \qquad\nabla\+k_{\}z_i}(z\);\la)\,\simeq\,D\+n_{z_i}(\la\);z)\,, &
D\+k_{\la_a}(z\);\la)\,\simeq\,\nabla\+n_{\]\la_a}(\la\);z)\,, &
\Tag{nD}
\\[8pt]
& \qquad\habla\+k_{\}z_i}(z\);\la)\,\simeq\,\Dh\+n_{z_i}(\la\);z)\,, &
\Dh\+k_{\la_a}(z\);\la)\,\simeq\,\habla\+n_{\]\la_a}(\la\);z)\,, &
\Tag{hD}
\\[8pt]
& Z\+k_{z_i}(z\);\la)\,\simeq\,N\+n_i(z)\>Q\+n_{z_i}(\la\);z)\,, \kern2em &
N\+k_a(\la)\>Q\+k_{\la_a}(z\);\la)\,\simeq\,Z\+n_{\la_a}(\la\);z)\,.
\kern-2em & \quad
\Tag{ZQ}
\cnn.2>
\end{alignat}
Here
\vvn-.2>
\begin{align}
N\+n_i(z)\, &{}=\]\prod_{1\le j<i}C\+n_{ji}(z_{ji}-\ka)
\}\prod_{i<j\le n}\bigl(C\+n_{ij}(z_{ij})\]\bigr)\vpb{-1}
\Tag{Nnz}
\\[-4pt]
\Text{and}
\nn.2>
N\+k_a(\la)\, &{}=\]\prod_{1\le b<a}C\+k_{ba}(\la_{ba}-\ka)
\}\prod_{a<b\le k}\bigl(C\+k_{ab}(\la_{ab})\]\bigr)\vpb{-1}.
\Tag{Nkl}
\end{align}
\end{theorem}
\vsk.1>
\nt
Equalities \(nD) and \(hD) for \difl/ operators are verified
in a straightforward way. Equalities \(ZQ) for \dif/ operators
follow from Theorem \ref{BR}\).
\begin{theorem}
\cite{TV4}
\label{BR}
For any \)$a\),b=1\lc k$, $a\ne b$, \>and any \>$i\),j=1\lc n$, $i\ne j$,
we have
\vvn.2>
$$
B\+k_{ab}(t)\>C\+k_{ab}(t)\,\simeq\,R\+n_{ab}(t)\,,\qqq
R\+k_{ij}(t)\,\simeq\,B\+n_{ij}(t)\>C\+n_{ij}(t)\,.
$$
\end{theorem}
\vsk.2>
Fix vectors $\lbi=(\)l_1\lc l_n)\in\Zp^{\)n}$ and
$\mbi=(\)m_1\lc m_k)\in\Zp^{\)k}$ \st/ $\sum_{i=1}^n\)l_i\)=\sum_{a=1}^k\)m_a$.
Let
\vvn-.2>
$$
\tsize\Zknlm\,=\,\bigl\lb\>
(d_{ai})_{\tsize\fratop{a\)=\)1\lc k}{\Lph ai\>=\)1\lc n}}\]\in\>\Zp^{\)kn}\;
\ \big|\ \sum_{a=1}^k\)d_{ai}=l_i\,,\quad\sum_{i=1}^n\)d_{ai}=m_a\)\bigr\rb\,.
$$
Denote by $\Plm\sub\PP_{kn}$ the span of all monomials
$\prod_{a=1}^k\)\prod_{i=1}^n\)x_{ai}^{\)d_{ai}}$ \st/ $(d_{ai})\in\Zknlm$.
\vvn.1>
Formulae \(VV)\), \(kiso)\,--\,\(right) and Proposition \ref{Pkn} imply that
$\Plm$ is isomorphic to each of the \wt/ subspaces
\vvn.4>
$$
\Vklm\quad\text{ and }\quad\Vnml\,.
\vv.4>
$$
The \iso/s are described in Proposition~\ref{Plm}.
\vsk.2>
Let $v_i\+k\],\,v_j\+n$ be \hwv/s of the respective modules
$V_{l_i}\+k\],\,V_{m_j}\+n$. For an indeterminate $y$ set $y\#0\}=\)1$
\,and \,$y\#s\}=\)y^s\}/\)s!$ \,for $s\in\Zpp$. For any \,$\db\in\Zknlm$ \,set
\,$x\#{\db\)}=\)\prod_{a=1}^k\)\prod_{i=1}^n\)x_{ai}\#{d_{ai}}\in\Plm$.
\begin{lemma}
\label{bases}
A basis of the \wt/ subspace $\Vklm$ is given by vectors
\begin{alignat}3
v_\db\+k\,={}&& \,\prod_{a=2}^k\bigl(e_{a1}\+k\bigr)\#{d_{a1}}v_1\+k\}
&{}\lox\prod_{a=2}^k\bigl(e_{a1}\+k\bigr)\#{d_{a\]n}}v_n\+k\,,
\kern1.5em && \db\)=\)(d_{ai})\in\Zknlm\,.\kern-1.2em
\Tag{vdk}
\\[6pt]
\Text{A basis of the \wt/ subspace $\Vnml$ is given by vectors}
\nn5>
v_\db\+n\,={}&& \,\prod_{i=2}^n\bigl(e_{i1}\+n\bigr)\#{d_{1i}}v_1\+n\}
&{}\lox\prod_{i=2}^n\)\bigl(e_{i1}\+n\bigr)\#{d_{ki}}v_k\+n\,,
&& \db\)=\)(d_{ai})\in\Zknlm\,.\kern-1.2em
\Tag{vdn}
\end{alignat}
\end{lemma}
\vsk.2>
\begin{proposition}
\label{Plm}
The \iso/s \(kiso) and \(niso) induce the \iso/s
\vvn.5>
\begin{alignat*}2
& \Vklm\,\to\,\Plm\,,\qqq && v_\db\+k\map\>x\#{\db\)},
\\[8pt]
& \Vnml\,\to\,\Plm\,,\qqq && v_\db\+n\map\>x\#{\db\)}.
\end{alignat*}
\end{proposition}
\vsk.4>
Since all \KZ/, \qKZ/, \DD/ and \qDD/ operators respect the \wt/ decomposition
of the corresponding tensor products of $\glk$ and \$\gln$-modules, they can
be restricted to \fn/s with values in \wt/ subspaces. Then one can read
Theorem~\ref{dualKZDD} as follows.
\vsk.1>
\begin{theorem}
\label{altKZDD}
Let $\pho$ be the \iso/ of \wt/ subspaces:
\vvn.4>
\begin{align}
\pho:\Vklm\,\to\,{}&\Vnml\,,
\Tag{pho}
\\[6pt]
\pho:v_\db\+k\to\,v_\db\+n\),\qqq & \)\db\)\in\Zknlm\,.\kern-2em
\notag
&\cnn.2>
\end{align}
Then for any \,$i=1\lc n$ \>and \>$a=1\lc k$ \>we have
\vvn.4>
\begin{alignat}2
\nabla\+k_{\}z_i}(z\);\la)\,&{}=\,\pho^{-1}\)D\+n_{z_i}(\la\);z)\,\pho\,,
&\qquad D\+k_{\la_a}(z\);\la)\,
&{}=\,\pho^{-1}\>\nabla\+n_{\]\la_a}(\la\);z)\,\pho\,,
\Tag{nDa}
\\[8pt]
\habla\+k_{\}z_i}(z\);\la)\,&{}=\,\pho^{-1}\)\Dh\+n_{z_i}(\la\);z)\,\pho\,,
& \Dh\+k_{\la_a}(z\);\la)\,
&{}=\,\pho^{-1}\>\habla\+n_{\]\la_a}(\la\);z)\,\pho\,,
\Tag{hDa}
\cnn.2>
\end{alignat}
\begin{align}
& Z\+k_{z_i}(z\);\la)\,=\,\pho^{-1}\)N\+n_i(z)\>Q\+n_{z_i}(\la\);z)\,\pho\,,
\Tag{ZQa}
\\[6pt]
& N\+k_a(\la)\>Q\+k_{\la_a}(z\);\la)\,=\,
\pho^{-1}\)Z\+n_{\la_a}(\la\);z)\,\pho\,.
\Tag{QZa}
\cnn.2>
\end{align}
Here $N\+n_i(z)$, $N\+k_a(\la)$ are given by formulae \(Nnz)\), \(Nkl)\).
\end{theorem}
Observe in addition that the restrictions of operators \(C) to the \wt/
subspaces are proportional to the identity operator:
\vvn.4>
\begin{gather}
C\+k_{ab}(t)\big|_{\Vklm}\>=\,\prod_{s=1}^{m_b}\>\frac{t-m_a\]+s-1}{t+s}\;,
\Tag{Ck}
\\[4pt]
C\+n_{ij}(t)\big|_{\Vnml}\>=\,\prod_{s=1}^{l_j}\>\frac{t-l_i\]+s-1}{t+s}\;,
\Tag{Cn}
\cnn.4>
\end{gather}
and are \raf/ of $t$.
\vsk.2>
Theorem~\ref{altKZDD} can be ``analytically continued'' \wrt/ $l_1, m_1$.
Namely, the theorem remains true if $l_1, m_1$ are complex numbers, while
all other numbers $l_2\lc l_n$, $m_2\lc m_k$ are still integers, and
${\sum_{i=1}^n\)l_i\)=\sum_{a=1}^k\)m_a}$. In this case
\vv.1>
the modules $V\+k_{l_1}$ and $V\+n_{m_1}$ are to be \irr/ \hwm/s with \hw/
${(l_1,0\lc 0)}$ and ${(m_1,0\lc 0)}$, \resp/, and the definition of $\Zknlm$
remains intact except that $d_{11}$ can be any number. Formulae \(vdk) and
\(vdn) make sense because they do not contain $d_{11}$, and Lemma~\ref{bases}
holds. Formulae \(Ck) and \(Cn) for ${a<b}$ and ${i<j}$ make sense for complex
$l_1, m_1$ as well, which is enough to obtain $N\+n_i(z)$ and $N\+k_a(\la)$ by
\(Nnz)\), \(Nkl)\). The ``\anco/'' of Theorem~\ref{altKZDD} will be useful in
application to identities of \hint/s of different dimensions.

\section{Hypergeometric \sol/s of the \KZv/ and dynamical \eq/s}
In the remaining part of the paper we will restrict ourselves to the case of
the Lie algebra $\glt$, which corresponds to $k=2$ in the previous sections.
\vsk.1>
Fix vectors $\lbi=(\)l_1\lc l_n)$ and $\mbi=(\)m_1,m_2)$
\vvn-.4>
\st/ $\sum_{i=1}^n\)l_i\)=\)m_1\]+m_2$ and $m_2\in\Zp$.
Let
\vvn.1>
$$
\tsize\Zlm\,=\,\bigl\lb\>
(d_1\lc d_n)\in\Zp^{\)n}\;\ \big|\ \sum_{i=1}^n\)d_i=m_2\,,\quad
d_i\le l_i\ \text{ if }\ l_i\]\in\Zp\)\bigr\rb\,.
$$
Given $d_1\lc d_n$, set \>$d_{\)<i}=\sum_{j=1}^{i-1}d_j$, \,$i=1\lc n$.
\vsk.1>
Consider the \wt/ subspace ${\Vlm}$. It has a basis given by vectors
\vvn.1>
$$
v_\db\,=\,e_{21}\#{d_1}v_{l_1}\]\lox\)e_{21}\#{d_n}v_{l_n}\,,\qqq
\db\)=\)(d_1\lc d_n)\in\Zlm\,,\kern-2em
\vvn.5>
$$
where $v_1\lc v_n$ are respective \hwv/s of the modules $V_{l_1}\lc V_{l_n}$.
\vsk.2>
Define the \em{master \fn/}
\begin{align*}
\Phi_r &(t_1\lc t_r;\zllb\))\,=\,
e^{\la_1\!\}\sum_{i=1}^n\}l_iz_i\)-\)(\la_1-\)\la_2)\sum_{a=1}^r\}t_a}\x{}
\\[4pt]
&{}\x\,
(\la_1\]-\la_2)^{-\)r}\!\]\prod_{1\le i<j\le n}\!\](z_i\]-z_j)^{l_il_j}
\,\prod_{a=1}^r\,\prod_{i=1}^n\,(t_a\]-z_i)^{-l_i}\!\]
\prod_{1\le a<b\le r}\!\](t_a\]-t_b)^2\,,\kern-1.6em
\cnn.4>
\end{align*}
and the \em{\wt/ \fn/}
$$
w_\db(t_1\lc t_r;\zn)\,=\,
\Sym\biggl[\ \prod_{i=1}^n\,\prod_{a=1}^{d_i}\,\frac1{t_{a\)+\)d_{<i}}\]-z_i}
\ \biggr]\,,\kern-1em
$$
where $\db=\)(d_1\lc d_n)\in\Zp^{\)n}$, \,$r\)=\sum_{i=1}^n\)d_i$, \,and
\vvn.2>
$$
\Sym f(t_1\lc t_r)\>=\sum_\si\,f(t_{\si_1}\lc t_{\si_r})\,.
$$
\vsk.3>
Fix a complex number $\ka$. Define a \$\Vlm\)$-valued \fn/
$U_\gm(\zlla)$ by the formula
\vvn.5>
\begin{alignat}2
\kern.6em U_\gm &{}(\zllm)\,={}
\Tag{Ugm}
\\[6pt]
&\>{}=\int_{\gm(\zlla)}\!\!\!\!
\bigl(\Phi_{m_2}(t_1\lc t_{m_2};\zllb\))\bigr)^{\]1/\ka}\x\}{}
\notag
\\[-1pt]
&& \Llap{{}\x\!\sum_{\db\in\Zlm}\!\!w_\db(t_1\lc t_{m_2};\zn)\,v_\db\;
d^{\>m_2}t\,} &.
\notag
\cnn.2>
\end{alignat}
The \fn/ depends on the choice of integration chains $\gm(\zlla)$. We assume
that for each $\zn$, $\la_1,\la_2$ the chain lies in $\C^{\>m_2}\}$ with
coordinates $t_1\lc t_{m_2}$, and the chains form a horizontal family of
\$m_2$-dimensional homology classes with respect to the multivalued \fn/
$\bigl(\Phi_{m_2}(t_1\lc t_{m_2};\zllb\))\bigr)^{\]1/\ka}\}$, see a more
precise statement below and in \cite{FMTV}\).
\begin{theorem}
\label{hgeom}
For any choice of the horizontal family $\gm$, the \fn/ $U_\gm(\zllm)$ is a
\sol/ of the \KZ/ and \DD/ \eq/s, see~\(KZ)\), \(DD)\), with values in $\Vlm$.
\end{theorem}
\vsk.1>
\nt
The theorem is a corollary of Theorem 3.1 in \cite{FMTV}\). For the \KZe/ at
$\la_1=\la_2$ the theorem follows from the results of \cite{SV}\), \cite{V}\).
\vsk.5>
There exist special horizontal families of integration chains in \(Ugm) labeled
by elements of $\Zlm$. They are described below. To simplify exposition we will
assume that $\Re\bigl((\la_1\]-\la_2)/\ka\bigr)>0$ and $\Im z_1\lsym<\Im z_n$.
\vsk.1>
Let $\db=\)(d_1\lc d_n)$, \,$r\)=\sum_{i=1}^n\)d_i$. Set
$\gm_\db(\zn)\,=\,\Cc_1\}\lx\Cc_r$, where $\Cc_1\lc\Cc_r$ is a collection of
non-intersecting oriented loops in $\C$ \st/ $d_i$ loops start at $+\)\infty$,
go around $z_i$, and return to $+\)\infty$, see the picture for $n=2$:
\vv1.1>\nl
\vbox{\begin{center}
\begin{picture}(222,85)
\put(40,15){\oval(10,10)[l]}
\put(40,15){\circle*{3}}
\put(40,10){\vector(1,0){152}}
\put(40,20){\line(1,0){152}}
\qbezier[16](46,5)(118,5)(190,5)
\qbezier[16](46,25)(118,25)(190,25)
\qbezier[4](30,15)(30,25)(40,25)
\qbezier[4](40,5)(30,5)(30,15)
\put(40,15){\oval(30,30)[l]}
\put(40,0){\vector(1,0){152}}
\put(40,30){\line(1,0){152}}

\put(31,65){\oval(10,10)[l]}
\put(31,65){\circle*{3}}
\put(31,60){\vector(1,0){161}}
\put(31,70){\line(1,0){161}}
\qbezier[17](37,55)(113,55)(190,55)
\qbezier[17](37,75)(113,75)(190,75)
\qbezier[4](21,65)(21,75)(31,75)
\qbezier[4](31,55)(21,55)(21,65)
\put(31,65){\oval(30,30)[l]}
\put(31,50){\vector(1,0){161}}
\put(31,80){\line(1,0){161}}

\put(13,5){$z_1$}
\put(4,55){$z_2$}
\put(197,80){$t_r$}
\put(197,70){$t_{d_1\]+1}$}
\put(197,30){$t_{d_1}$}
\put(197,20){$t_1$}
\end{picture}
\vsk>
Picture 1. The contour $\gm_\db$.
\end{center}}
\vsk.9>
One can see that for any ${\db\in\Zlm}$ the family of chains $\gm_\db$ is
horizontal. Therefore, the \fn/ $\Udzlm=U_{\gm_\db}(\zllm)$ is a \sol/ of
systems \(KZ) and \(DD)\). A univalued branch of the integrand in \(Ugm)
is fixed by assuming that at the point of $\gm_\db$ where all numbers
$t_{a\)+\)d_{<i}}\]-z_i$, $i=1\lc n$, $a=1\lc d_i$, are negative one has
\vvn.4>
$$
-\)\pi<\arg\)(z_i\]-z_j)<0\,,\ \quad
-\)2\pi<\arg\)(t_a\]-z_i)<0\,,\ \quad -\)\pi<\arg\)(t_a\]-t_b)\le 0\,,
\vv.4>
$$
for $i=1\lc n$, $j=i+1\lc n$, \>$a=1\lc m_2$, $b=a+1\lc m_2$.
\vsk.1>
The \sol/ $\Udzlm$ is distinguished by the following property.
\begin{theorem}
\label{Udas}
Let $\Im(z_i-z_{i+1})\to\)-\)\infty$ for all $i=1\lc n-1$.
Then for any $\db\in\Zlm$ one has
\vvn.3>
\begin{align}
U_\db &(\zllm)\,={}
\Tag{Udas}
\\[7pt]
&\){}=\,(2\pi i)^{m_2}\,e^{\pi i\)\xi_\db(\lbi)/\ka}\,
\bigl(\)\Xi_{\)\db}(\zllm)\bigr)^{\]1/\ka}\)\x{}
\notag
\\[5pt]
&\hp{{}={}}\x\,\prod_{j=1}^n\,\prod_{s=0}^{d_j-1}
\frac{\Gm(-\)1/\ka)}{\Gmb{1+(l_j\]-s)/\ka}\,\Gmb{-\)(s+1)/\ka}}
\ \bigl(v_\db+o(1)\bigr)\kern-1em
\notag
\cnn.5>
\end{align}
where \,$\xi_\db(\)\lbi\))\>=\!\}\sum_{1\le i\le j\le n}\!\]l_i\)d_j$ \,and
\vvn-.4>
\begin{align*}
\quad\Xi_{\)\db}& (\zllm)\,=\,\ka^{-\)m_2}\,
e^{\la_1\!\]\sum_{i=1}^n\}z_i(l_i-d_i)\)+\)\la_2\!\]\sum_{i=1}^n\}z_id_i}\)\x{}
\\[5pt]
& {}\>\x\,\bigl((\la_1\]-\la_2)/\ka\bigr)^{\)\sum_{i=1}^n\!d_i(l_i-\)d_i)}
\!\!\prod_{1\le i<j\le n}\!(z_i\]-z_j)^{(l_i-\)d_i)(l_j-\)d_j)\)+\)d_id_j}\>.
\cnn-.4>
\end{align*}
\end{theorem}
\vsk.6>
Theorem \ref{Udas} implies that the set of \sol/s $U_\db$, $\db\in\Zlm$, of
systems \(KZ) and \(DD) is complete, that is, any \sol/ of those systems taking
values in $\Vlm$ is a linear combination of \fn/s $U_\db$.
\vsk.1>
There is a similar statement for \as/s of $\Udzlm$ \wrt/ $\la_1, \la_2$.
\begin{theorem}
\label{Udla}
Let $\Re\bigl((\la_1\]-\la_2)/\ka\bigr)\to\)+\)\infty$.
Then for any $\db\in\Zlm$ formula \(Udas) holds.
\end{theorem}
The proof of Theorems~\ref{Udas} and \ref{Udla} uses the following
Selberg-type integral
\vvn.4>
\begin{align}
\int_{\gm_m\!}{} & \;e^{-\)\nu\!\sum_{a=1}^m\!s_a}\)
\prod_{a=1}^m\>(-\)s_a)^{-1-\)l/\ka}\!\!
\prod_{1\le a<b\le m}\!(s_a\]-s_b)^{2/\ka}\,d^{\>m}\}s\,={}
\Tag{selb}
\\[1pt]
& {}=\,(-\)2\pi i)^m\>\nu^{m\)(l-\)m\)+\)1)/\ka}\;\prod_{j=0}^{m-1}\,
\frac{\Gm(1-1/\ka)}{\Gmb{1+(l-j)/\ka}\,\Gmb{1-(j+1)/\ka}}\;,\kern-1em
\notag
\cnn.2>
\end{align}
where $\Re\nu>0$,
\>${\gm_m=\)\bigl\lb(s_1\lc s_m)\in\C^{\>m}\ |\ s_a\in\Cc_a}$,
$a=1\lc m\)\bigr\rb$, \>and $\Cc_1\lc \Cc_m$ are non-intersecting oriented
loops in $\C$ which start at $+\)\infty$, go around zero, and return to
$+\)\infty$, the loop $\Cc_a$ being inside $\Cc_b$ for $a<b$,
see the picture:\vvn.9>\nl
\vbox{\begin{center}
\begin{picture}(224,35)
\put(40,15){\oval(10,10)[l]}
\put(40,15){\circle*{3}}
\put(40,10){\vector(1,0){152}}
\put(40,20){\line(1,0){152}}
\qbezier[16](46,5)(118,5)(190,5)
\qbezier[16](46,25)(118,25)(190,25)
\qbezier[5](30,15)(30,25)(40,25)
\qbezier[4](40,5)(30,5)(30,15)
\put(40,15){\oval(30,30)[l]}
\put(40,0){\vector(1,0){152}}
\put(40,30){\line(1,0){152}}
\put(15,5){$0$}
\put(197,30){$s_m$}
\put(197,20){$s_1$}
\end{picture}
\vsk>
Picture~2. The contour $\gm_m$.
\end{center}\vsk.7>}
A univalued branch of the integrand in \(selb) is fixed by assuming that
at the point of $\gm_m$ where all numbers $s_1\lc s_m$ are negative one
has $\arg\)(-\)s_1)\lsym=\arg\)(-\)s_m)=0$ and $\arg\)(s_a\]-s_b)=0$ for
$1\le a<b\le m$.
\vsk.2>
The construction of \hgeom/ \sol/s of the \tri/ \KZe/s \(KZt) and
the \dif/ dynamical \eq/s \(qDD) is similar. We describe it below.
\vsk.1>
Define the master \fn/
\vvn.3>
\begin{align}
\label{mastt}
\Psi_r &(t_1\lc t_r;\zllm)\,={}
\\[5pt]
&{}=\,\prod_{i=1}^n\,z_i^{\)l_i(\la_1-\)m_1+\)l_i/2)}\!\}
\prod_{1\le i<j\le n}\!\](z_i\]-z_j)^{l_il_j}\x{}
\notag
\\[4pt]
& {}\>\x\,\prod_{a=1}^r\,t_a^{\)\la_2-\la_1+\)m_1-\)m_2\)+1}
\,\prod_{a=1}^r\,\prod_{i=1}^n\,(t_a\]-z_i)^{-l_i}\!\]
\prod_{1\le a<b\le r}\!\](t_a\]-t_b)^2\,,\kern-2em
\notag
\cnn-.6>
\end{align}
\goodbreak\vsk1.2>\nt
and a \$\Vlm\)$-valued \fn/
\vvn.5>
\begin{alignat}2
\kern.6em \Ut_\dl &{}(\zllm)\,={}
\Tag{Udl}
\\[4pt]
&\>{}=\int_{\dl(\zlla)}\!\!\!\!
\bigl(\Psi_{m_2}(t_1\lc t_{m_2};\zllm)\bigr)^{\]1/\ka}\x\}{}
\notag
\\[-4pt]
&& \Llap{{}\x\!\sum_{\db\in\Zlm}\!\!w_\db(t_1\lc t_{m_2};\zn)\,v_\db\;
d^{\>m_2}t\,} &.
\notag
\cnn-.5>
\cnn.7>
\end{alignat}
The \fn/ depends on the choice of integration chains $\dl(\zlla)$. We assume
that for each $\zn$, $\la_1,\la_2$ the chain lies in $\C^{\>m_2}\}$ with
coordinates $t_1\lc t_{m_2}$, and the chains form a horizontal family of
\$m_2$-dimensional homology classes with respect to the multivalued \fn/
$\bigl(\Psi_{m_2}(t_1\lc t_{m_2};\zllm)\bigr)^{\]1/\ka}\}$, see a more
precise statement below and in \cite{MV}\).
\begin{theorem}
\label{hgeomt}
For any choice of the horizontal family $\dl$, the \fn/ $U_\dl(\zllm)$ is
a \sol/ of the \tri/ \KZ/ and \rat/ \qDD/ \eq/s, see~\(KZt)\), \(qDD)\),
with values in $\Vlm$.
\end{theorem}
\nt
The theorem is a direct corollary of results in \cite{MV}\). Another way
of writing down \hgeom/ \sol/s of the \tri/ \KZe/s is given in \cite{Ch1}\).
\vsk.5>
There exist special horizontal families of integration chains in \(Udl) labeled
by elements of $\Zlm$. They are described below. To simplify exposition we will
assume that $\Re\bigl((\la_1\]-\la_2)/\ka\bigr)$ is large positive and
$\arg\)z_1\lsym<\arg\)z_n<{\arg\)z_1\]+2\pi}$, that is, all the ratios
$z_i/z_j$ for $i\ne j$ are not real positive, and $\zn$ are ordered
counterclockwise. Recall that all $\zn$ are nonzero.
\vsk.1>
Let $\db=\)(d_1\lc d_n)$, \,$r\)=\sum_{i=1}^n\)d_i$. Set
$\dl_\db(\zn)\,=\,\Cc_1\}\lx\Cc_r$, where $\Cc_1\lc\Cc_r$ is a collection of
non-intersecting oriented loops in $\C$ \st/ $d_i$ loops start at infinity in
the direction of $z_i$, go around $z_i$, and return to infinity in the same
direction, see the picture for $n=2$:
\nl
\vbox{\begin{center}
\begin{picture}(152,135)
\put(45,15){\oval(10,10)[l]}
\put(45,15){\circle*{3}}
\put(45,10){\vector(1,0){120}}
\put(45,20){\line(1,0){120}}
\qbezier[16](51,5)(105,5)(163,5)
\qbezier[16](51,25)(105,25)(163,25)
\qbezier[4](35,15)(35,25)(45,25)
\qbezier[4](45,5)(35,5)(35,15)
\put(45,15){\oval(30,30)[l]}
\put(45,0){\vector(1,0){120}}
\put(45,30){\line(1,0){120}}

\put(0,35){\circle*{3}}
\qbezier(-4,31)(-8,35)(-4.2,39)
\qbezier(-4,31)(0,28.1)(4,31.2)
\put(-4,39){\line(1,1){69}}
\put(4,31){\vector(1,1){69}}
\qbezier(-11,24)(-22,35)(-11.2,46)
\qbezier(-11,24)(0,14.4)(11,24.3)
\put(-11,46){\line(1,1){69}}
\put(11,24){\vector(1,1){69}}
\qbezier[10](11,31)(43,63)(75,95)
\qbezier[10](-4,46)(28,78)(60,110)
\qbezier[4](7,27)(0,21)(-8,27)
\qbezier[4](-8,42)(-15,35)(-8,27)

\put(16,13){$z_1$}
\put(-31,33){$z_2$}
\put(59,120){$t_r$}
\put(68,111){$t_{d_1\]+1}$}
\put(170,30){$t_{d_1}$}
\put(170,20){$t_1$}

\put(-35,13){$0$}
\put(-20,15){\circle*{3}}
\end{picture}
\vsk>
Picture~3. The contour $\dl_\db$.
\end{center}\vsk-.4>}
\vsk1.4>
One can see that for any ${\db\in\Zlm}$ the family of chains $\dl_\db$ is
horizontal. Therefore, the \fn/ $\Utdzlm=\Ut_{\dl_\db}(\zllm)$ is a \sol/
of systems \(KZt) and \(qDD)\). A univalued branch of the integrand in \(Udl)
is fixed by assuming that at the point of $\dl_\db$ where all ratios
$t_{a\)+\)d_{<i}}/z_i$, $i=1\lc n$, $a=1\lc d_i$, are real and belong to
$(0\>,1)$ one has
\vvn.2>
$$
\arg\)t_{a\)+\)d_{<i}}\]=\>\arg\)z_i\,,\qquad i=1\lc n,\quad a=1\lc d_i\,,
\kern-3em
\vv-.4>
$$
and
\vvn-1.2>
\begin{gather*}
-\)\pi<\arg\)(z_i\]-z_j)-\)\arg\)z_i<\pi\,,
\\[4pt]
-\)2\pi<\arg\)(t_a\]-z_i)-\)\arg\)z_i<0\,,
\\[4pt]
-\)\pi<\arg\)(t_a\]-t_b)-\)\arg\)t_a<\pi\,,
\cnn.4>
\end{gather*}
for $i=1\lc n$, $j=i+1\lc n$, \>$a=1\lc m_2$, $b=a+1\lc m_2$.
Recall, it is assumed that $\arg\)z_1\lsym<\arg\)z_n<\arg\)z_1\]+2\pi$,
\vsk.2>
There is an analogue of Theorem~\ref{Udas} which describes \as/s of the \fn/s
$\Utdzlm$ as $z_i/z_{i+1}\to 0$ for all $i=1\lc n-1$. The corresponding
formulae are similar to \(Udas), but more involved. Asymptotics of $\Ut_\db$
\wrt/ $\la_1, \la_2$ are as follows.
\begin{theorem}
\label{Utas}
Let $(\la_1\]-\la_2)/\ka\to\)+\)\infty$. Then
\vvn.3>
\begin{align}
\Ut_\db &(\zllm)\,={}
\Tag{Utas}
\\[7pt]
&\){}=\,(2\pi i)^{m_2}\,e^{\pi i\)\xi_\db(\lbi)/\ka}\,
\bigl(\)\Xt_{\)\db}(\zllm)\bigr)^{\]1/\ka}\)\x{}
\notag
\\[5pt]
&\hp{{}={}}\x\,\prod_{j=1}^n\,\prod_{s=0}^{d_j-1}
\frac{\Gm(-\)1/\ka)}{\Gmb{1+(l_j\]-s)/\ka}\,\Gmb{-\)(s+1)/\ka}}
\ \bigl(v_\db+o(1)\bigr)\kern-.6em
\notag
\cnn.3>
\end{align}
where \,$\xi_\db(\)\lbi\))\>=\!\}\sum_{1\le i\le j\le n}\!\]l_i\)d_j$ \,and
\vvn-.4>
\begin{align*}
& \Xt_{\)\db}(\zllm)\,=\,
\bigl((\la_1\]-\la_2)/\ka\bigr)^{\)\sum_{i=1}^n\!d_i(l_i-\)d_i+\)1)}\x{}
\\[5pt]
& {}\}\x\,\prod_{i=1}^n z_i^{(\la_1-\)m_1)(l_i-\)d_i)\)+\)(\la_2-\)m_2)\)d_i+
\)((l_i-\)d_i)^2+\)d_i^2)/2}
\!\!\prod_{1\le i<j\le n}\!(z_i\]-z_j)^{(l_i-\)d_i)(l_j-\)d_j)\)+\)d_id_j}\>.
\end{align*}
\end{theorem}
\vsk.4>
The construction of \hgeom/ \sol/s of the \qKZe/s \(KZt) and \tri/ \DD/ \eq/s
\(DDt) goes along the same lines as for \hgeom/ \sol/s of the \KZ/ and \rat/
dynamical \eq/s, but instead of ordinary \hint/s it employs \q-\hint/s of
Mellin\>-\)Barnes type, see~\cite{TV2}\).
\vsk>
\goodbreak
\vsk-.9>
Define the \em{\q-master \fn/}
\vvn.4>
\begin{align}
\label{qmaster}
\kern.6em \Phh_r & (t_1\lc t_r;\zllb\);\ka)\,={}
\\[4pt]
& {}=\,\la_1^{(r\]/2\,+\}\sum_{i=1}^n\}(z_il_i-\)l_i^2/2)\)-\}
\sum_{a=1}^r\}t_a)/\ka}\la_2^{(r\]/2\,+\}\sum_{a=1}^r\}t_a)/\ka}
(\la_1\]-\la_2)^{-\)r\]/\]\ka}\>\x{}
\notag
\\[6pt]
& \hp{{}={}}\}\x\,\prod_{a=1}^r\,
\prod_{i=1}^n\,\frac{\Gmb{(t_a\]-z_i)/\ka}}{\Gmb{(t_a\]-z_i\]+l_i)/\ka}}
\prod_{1\le a<b\le r}
\frac{\Gmb{(t_a\]-t_b\]+1)/\ka}}{\Gmb{(t_a\]-t_b\]-1)/\ka}}\;,
\notag
\cnn.4>
\end{align}
the \em{\rat/ \wt/ \fn/}
\vvn.3>
\begin{align*}
& \wh_\db(t_1\lc t_r;\zn;\lbi\))\,=
\prod_{1\le a<b\le r}\frac{t_a\]-t_b}{t_a\]-t_b\]-1}\;\x{}
\\
&{}\}\x\,\Sym\biggl[\ \prod_{j=1}^n\,\prod_{a=1}^{d_j}\,
\Bigl(\>\frac1{t_{a\)+\)d_{<j}}\]-z_j\]+l_j}\>
\prod_{1\le\)p<j}\frac{t_{a\)+\)d_{<j}}\]-z_p}{t_{a\)+\)d_{<j}}\]-z_p\]+l_p}\>
\Bigr)\!\}\prod_{1\le a<b\le r}\!\!\frac{t_a\]-t_b\]-1}{t_a\]-t_b}\ \biggr]\,,
\kern-.4em
\cnn.3>
\end{align*}
where $\db=\)(d_1\lc d_n)\in\Zp^{\)n}$, \,$r\)=\sum_{j=1}^n\)d_j$, \,and
\vvgood
the \em{\tri/ \wt/ \fn/}
\begin{alignat*}2
& \kern-1em W_\db(t_1\lc t_r;\zn;\lbi\))\,=
\prod_{1\le a<b\le r}\frac{\sink{t_a\]-t_b}}{\sink{t_a\]-t_b\]-1}}\;\x{}
\\[4pt]
& \kern-1em {}\}\!\x\,\Sym\biggl[\ \prod_{j=1}^n\,\prod_{a=1}^{d_j}\,
\Bigl(\>\frac{e^{\)\pi\)i\)(z_j-\)t_{a\)+\)d_{<j}})/\ka}}
{\sink{t_{a\)+\)d_{<j}}\]-z_j\]+l_j}}\>\prod_{1\le\)p<j}
\frac{\sink{t_{a\)+\)d_{<j}}\]-z_p}}{\sink{t_{a\)+\)d_{<j}}\]-z_p\]+l_p}}
\>\Bigr)\>\x\!\}{}\kern-1.6em &&
\\[3pt]
&&& \Llap{{}\x\!
\prod_{1\le a<b\le r}\!\!\frac{\sink{t_a\]-t_b\]-1}}{\sink{t_a\]-t_b}}
\ \biggr]\,.\kern-1.6em}
\cnn.1>
\end{alignat*}
\vsk.3>
For simplicity of exposition from now on we assume that $\ka$ is a real
positive number and the ratio $\la_2/\la_1$ is not real positive.
For any ${\db\in\Zlm}$ define a \$\Vlm\)$-valued \fn/
\vvn.4>
\begin{alignat}2
& \Uh_\db(\zllm)\,={}
\Tag{Uhdb}
\\[4pt]
{}=\int_{\II\)(\zn;\)\lbi)}\!\!\!\!{}
& \Phh_{m_2}(t_1\lc t_{m_2};\zllb\))\,
W_\db(t_1\lc t_{m_2};\zn;\lbi\))\,\x{}\!\}&
\notag
\\[-4pt]
&&& \Llap{{}\x\!\sum_{\pb\in\Zlm}\!\!\wh_\pb(t_1\lc t_{m_2};\zn;\lbi)
\,v_\pb\,\>d^{\>m_2}t\,,}
\notag
\cnn.2>
\end{alignat}
the integration contour $\II\)(\zn;\lbi)$ being described below.
For the factors $(\la_2/\la_1)^{\)t_a\]/\ka}$ in the integrand
it is assumed that $0<\arg\)(\la_2/\la_1)<2\pi$.
\vsk.1>
The integral in \(Uhdb) is defined by \anco/ \wrt/ $\zn$ and
$\lbi=(l_1\lc l_n)$ from the region where $\Re z_1\lsym=\Re z_n\]=0$
and $\Re l_i<0$ for all $i=1\lc n$. In that case
\vvn.2>
$$
\II\)(\zn;\lbi\))\)=\>\bigl\lb(t_1\lc t_{m_2})\in\C^{\>m_2}
\ |\ \Re t_1\lsym=\Re t_{m_2}\]=\eps\>\bigr\rb
\vv.3>
$$
where $\eps$ is a positive number less then
$\min\)(-\]\Re l_1\)\lc\]-\]\Re l_n)$. In the considered region of parameters
the integrand is well defined on $\II\)(\zn;\lbi\))$ and the integral is
convergent. It is also known that $\Uh_\db(\zllm)$ can be analytically
continued to a value of $\lbi$ in $\Zp^{\>n}$ and generic values of $\zn$,
if $\db\in\Zlm$ at that point, and the \anco/ is given by the integral over
a suitable deformation of the imaginary plane
$\bigl\lb(t_1\lc t_{m_2})\in\C^{\>m_2}
\ |\ \Re t_1\lsym=\Re t_{m_2}\]=0\>\bigr\rb$, see \cite{MuV}\).

\begin{theorem}
\label{qhgeom}
For any $\db\in\Zlm$ the \fn/ $\Uh_\db(\zllm)$ is a \sol/ of the \rat/ \qKZ/
and \tri/ \DD/ \eq/s, see~\(qKZ)\), \(DDt)\), with values in $\Vlm$.
\end{theorem}
\nt
The part of the theorem concerning the \qKZe/s is a direct corollary of
the construction of \q-\hgeom/ \sol/s of the \qKZe/s given in \cite{TV2}\),
\cite{MuV}\). The part of the theorem on the \tri/ \DD/ \eq/s is obtained in
\cite{TV8}\).
\vsk.5>
The \sol/ $\Uh_\db(\zllm)$ of systems \(qKZ) and \(DDt) is distinguished
by the following property.
\begin{theorem}
\label{Uhas}
Let $\Re(z_i-z_{i+1})\to\)+\)\infty$ for all $i=1\lc n-1$.
Then for any $\db\in\Zlm$ one has
\vvn.3>
\begin{align}
\Uh_\db &(\zllm)\,={}
\Tag{Uhas}
\\[7pt]
&\){}=\,(-2i)^{m_2}\,m_2!\;e^{\pi i\)\zt_\db(\lbi)/\ka}\,
\bigl(\)\Xh_{\)\db}(\zllm)\bigr)^{\]1/\ka}\)\x{}
\notag
\\[5pt]
&\hp{{}={}}\x\,\prod_{j=1}^n\)\biggl(d_j!\,\prod_{s=0}^{d_j-1}\>
\frac{\Gmb{(s-l_j)/\ka}\,\Gmb{1+(s+1)/\ka}}{\Gm(1+1/\ka)}\>\biggr)
\bigl(v_\db+o(1)\bigr)\,.
\notag
\cnn.3>
\end{align}
where \,$\zt_\db(\)\lbi\))\)=\sum_{i=1}^n d_j(2\)l_j\]-d_j\]+1)/2$,
\,and
\vvn-.2>
\begin{align*}
\Xh_{\)\db}(\zllm)\,=\,
\la_1^{\)\sum_{i=1}^n\!(z_i(l_i-d_i)-\)l_i^2/2\>+\)d_i^2/2)}
\la_2^{\)\sum_{i=1}^n\!d_i(z_i-\)l_i+\)d_i/2)}\x{}\!\}&
\\[5pt]
{}\x\,(\la_1\]-\la_2)^{\)\sum_{i=1}^n\!d_i(l_i-\)d_i)}\!\!
\prod_{1\le i<j\le n}\!
\bigl((z_i\]-z_j)/\ka\bigr)^{(l_i-\)d_i)(l_j-\)d_j)\)+\)d_id_j-\)l_il_j} &
\cnn.2>
\end{align*}
Recall that $\ka$ is assumed to be a real positive number.
\end{theorem}
Theorem \ref{Udas} implies that the set of \sol/s $\Uh_\db$, $\db\in\Zlm$, of
systems \(qKZ) and \(DDt) is complete, that is, any \sol/ of those systems
taking values in $\Vlm$ is a linear combination of \fn/s $U_\db$.
\vsk.2>
The proof of Theorem~\ref{Udas} uses the following Selberg-type integral
\vvn.3>
\begin{align}
\!\!\int_{\II_m(l)\!}(-x)^{\>\sum_{a=1}^m\!s_a}\)
\prod_{a=1}^m\,\Gm(s_a)\,\Gm(\)-\)s_a\]-l/\ka)\)
\prod_{\fratop{a,b=1}{a\ne b}}^m\]
\frac{\Gm(s_a\]-s_b\]+1/\ka)}{\Gm(s_a\]-s_b)}\ d^{\>m}\}s\,={}&
\Tag{qselb}
\\[6pt]
{}=\,(2\)\pi i)^m\,(-\)x)^{\)(m-1-2\)l)\)m/2\ka}\,
(1-x)^{\)m\)(l-m+1)/\ka}\,\x{} &
\notag
\\[3pt]
{}\x\,\prod_{j=0}^{m-1}\,
\frac{\Gmb{(j-l)/\ka}\,\Gmb{1+(j+1)/\ka}}{\Gm(1+1/\ka)}\;, &
\notag
\cnn.3>
\end{align}
where ${-\pi<\arg\)(-\)x)<\pi}$ and ${-\)\pi<\arg\)(1-x)<\pi}$. The integral
is defined by \anco/ from the region where $\ka$ is real positive and $\Re l$
is negative. In that case
\vvn.4>
$$
\II_m(l)\>=\,\bigl\lb(s_1\lc s_m)\in\C^{\>m}\ |\ \,
\Re s_1\lsym=\)\Re s_m=\)-\)\Re l/2\)\bigr\rb\,.
\vv.3>
$$
In the considered region of parameters the integrand in \(qselb) is well
defined on $\II_m(l)$ and the integral is convergent, see \cite{TV1}\).

\section{Duality for \hgeom/ and \boldmath\q-\hint/s}
In this section we consider the $\glkn$ duality for the case of $k=n=2$, and
apply the results of the previous sections to obtain identities for \hgeom/
and \q-\hint/s of different dimensions.
\vsk.1>
Further on we fix complex numbers $l_1, m_1$ and nonnegative integers
$l_2, m_2$ \st/ $l_1\]+l_2=\)m_1\]+m_2$. Set $\lbi=(l_1,l_2)$ and
$\mbi=(m_1,m_2)$.
\vsk.1>
Let $V_l$ be the \irr/ \hw/ \$\glt$-module with \hw/ $(l,0)$ and \hwv/ $v_l$.
The \wt/ subspace $\Vlmt$ has a basis given by vectors
\vvn.4>
$$
v_b(\lmb)\,=\,\)\frac1{(m_2\]-b)!\;b\)!}\;
e_{21}^{(m_2-b)}v_{l_1}\}\ox e_{21}^bv_{l_2}\,,\qquad b=0\lc\min\)(l_2,m_2)\,,
\vv.4>
$$
provided that $l_1$ is not a nonnegative integer or $m_2\le l_1$. Otherwise,
the vectors $v_0(\lmb)\lc v_{m_2-\)l_1-\)1}(\lmb)$ equal zero and the basis is
given by the rest of the vectors $v_b(\lmb)$. Say that $b$ is admissible if
$v_b(\lmb)\ne 0$.
\vsk.1>
The \wt/ subspaces $\Vlmt$ and $\Vmlt$ are isomorphic.
The \iso/ $\pho:\Vlmt\to\Vmlt$, \cf. \(pho)\), sends the vector $v_b(\lmb)$
to $v_b(\mlb\))$.
\vsk.1>
Given an admissible integer $b$ let ${\db=(m_2\]-b,b)}$ and
${\db'=(l_2\]-b,b)}$. Consider \$\Vlmt\)$-valued \fn/s
\vvn.4>
\begin{align*}
& U_b(\zlt;\lmb)\,=\,U_\db(\zlt;\lmb)
\\[-5pt]
\Text{and}
\nn-5>
& U'_b(\lzt;\mlb)\,=\,\pho^{-1}\bigl(U_{\db'}(\lzt;\mlb\))\bigr)\,,
\cnn.3>
\end{align*}
where the \fn/s $U_\db$ and $U_{\db'}$ are defined in Section~7, \cf. \(Ugm)
and Picture~1.
\begin{theorem}
\cite{TV6}
\label{dualhint}
For any $b=0\lc\min\)(l_2,m_2)$ one has
\vvn.4>
\begin{gather}
A_b(\lmb)\>U_b(\zlt;\lmb)\,=\,A_b(\mlb)\>U'_b(\lzt;\mlb)
\Tag{dualhint}
\cnn.1>
\end{gather}
where
\vvn-.2>
\begin{align*}
A_b(\lmb)\, &{}=\,(-2\)i)^{-m_2}\>\ka^{(m_1+\)1)m_2/\ka}\,
e^{-\pi i\)(m_1+\)m_2-\)b)m_2/\ka}\>\x{}
\notag
\\[3pt]
& \>{}\x\prod_{s=0}^{m_2-b-1}\!\frac1{\sink{s+1}}\,\,
\prod_{s=0}^{m_2-1}\>
\frac{\Gmb{1+(l_1\]-s)/\ka}}{\Gm(-\)1/\ka)\,\Gmb{1+(s+1)/\ka}}\;.\kern-1em
\notag
\end{align*}
\end{theorem}
\vsk.4>
\nt
The idea of the proof of the statement is as follows. By Theorems~\ref{hgeom}
and \ref{altKZDD} the \fn/s $U_b(\zlt;\lmb)$ and $U'_b(\lzt;\mlb)$ are \sol/s
of the \rat/ \difl/ \KZ/ and dynamical \eq/s \(KZ) and \(DD)\).
Theorem~\ref{Udas} implies that the \fn/s $U_b(\zlt;\lmb)$ with admissible
$b\}$'s form a complete set of \sol/s, which means that the \fn/s
$U'_b(\lzt;\mlb)$ are their linear combinations. The transition coefficients
can be found by comparing \as/s of $U_b$ and $U'_b$ as $z_1\]-z_2$ goes to
infinity, see Theorems~\ref{Udas} and \ref{Udla}.
\begin{remark}
Equality \(dualhint) holds for vector-valued \fn/s. That is, it contains
several identities of the form: a \hint/ of dimension $m_2$ (a coordinate
of $U_b$) equals a \hint/ of dimension $l_2$ (the corresponding coordinate
of $U'_b$).
\end{remark}
Consider \$\Vlmt\)$-valued \fn/s
\vvn.4>
$$
\Uh_b(\zlt;\lmb)\,=\,\Uh_\db(\zlt;\lmb)
\vvn-.5>
$$
and
\begin{align*}
\Ut'_b & (\lzt;\mlb)\,={}
\\[3pt]
&{}=\,\prod_{s=0}^{l_2-1}\,
\frac{\Gmb{(z_1\]-z_2+s-l_1)/\ka}}{\Gmb{(z_1\]-z_2+s+1)/\ka}}
\ \;\pho^{-1}\bigl(\Ut_{\db'}(\lzt;\mlb\))\bigr)\,,\kern-1em
\cnn.2>
\end{align*}
where ${\db=(m_2\]-b,b)}$, \,${\db'=(l_2\]-b,b)}$, \,and the \fn/s $\Uh_\db$,
$\Ut_{\db'}$ are defined in Section~7, \cf. \(Uhdb)\), \(Udl) and Picture~3.
\begin{theorem}
\cite{TV6}
\label{dualqhint}
For any $b=0\lc\min\)(l_2,m_2)$ one has
\vvn.4>
\begin{gather}
\Ah_b(\lmb)\>\Uh_b(\zlt;\lmb)\,=\,\At_b(\mlb)\>\Ut'_b(\lzt;\mlb)
\Tag{dualqhint}
\cnn.1>
\end{gather}
where
\vvn-.1>
$$
\Ah_b(\lmb)\,=\,(2\pi i)^{-m_2}
\prod_{s=0}^{m_2-b-1}\!\sink{l_1\]-s}\,
\prod_{s=0}^{m_2-1}\>
\frac{\Gmb{1+(l_1\]-s)/\ka}\,\Gm(1+1/\ka)}{\Gmb{1+(s+1)/\ka}}\;,
\vv-.2>
$$
and
\begin{align*}
\At_b(\lmb)\, &{}=\,(2\pi i)^{-l_2}\,e^{\pi i\)(-b^2+
\)b\)(l_2-\)l_1)\)+\)l_1m_2\)-\)l_2m_1-\)m_2(m_2\)-\)1)/2)/\ka)}\>\x{}
\notag
\\[3pt]
& \>{}\x\,\prod_{s=0}^{b-1}\,\frac1{\sink{s+1}}\,\,
\prod_{s=0}^{l_2-1}\>
\frac{\Gmb{1+(m_1\]-s)/\ka}\,\Gmb{-\)(s+1)/\ka}}{\Gm(-\)1/\ka)}\;.\kern-1em
\notag
\end{align*}
\end{theorem}
\vsk.4>
\nt
The idea of the proof is similar to that of Theorem~\ref{dualhint}.
By Theorems~\ref{qhgeom} and \ref{hgeomt}, \ref{altKZDD} the \fn/s
$\Uh_b(\zlt;\lmb)$ and $\Ut'_b(\lzt;\mlb)$ are \sol/s of the \rat/ \qKZe/s
\(qKZ)\). \)Theorem~\ref{Uhas} implies that the \fn/s
$U_b(\zlt;\lmb)$ with admissible $b\}$'s form a complete set of \sol/s over
the field of \$\ka\)$-periodic \fn/s of $z_1, z_2$ ($\la_1, \la_2$ are treated
as parameters in the present consideration). Therefore, the \fn/s
$U'_b(\lzt;\mlb)$ as \fn/s of $z_1, z_2$ are linear combinations of
$U_b(\zlt;\lmb)$ with periodic coefficients. The coefficients can be found
by comparing \as/s of $U_b$ and $U'_b$ as $z_1\]-z_2$ goes to infinity,
see Theorems~\ref{Uhas} and \ref{Utas}.
\par
Though one does not need the fact that the \fn/s $\Uh_b(\zlt;\lmb)$ and
$\Ut'_b(\lzt;\mlb)$ solve the \tri/ \DD/ \eq/s \(DDt) in the proof of
Theorem~\ref{dualqhint}, this fact is reflected in formula \(dualqhint)
--- the coefficients $\Ah_b(\lmb)$ and $\At_b(\lmb)$ do not depend on
$\la_1, \la_2$.
\begin{remark}
Similar to \(dualhint), equality \(dualqhint) contains several identities of
the form: a \q-\hint/ of dimension $m_2$ (a coordinate of $\Uh_b$) equals
a \hint/ of dimension $l_2$ (the corresponding coordinate of $\Ut'_b$).
\end{remark}
For $l_2=m_2=1$ formula \(dualqhint) yields the classical equality of integral
\rep/s of the Gauss \hgeom/ \fn/ $_2F_1$. For instance, taking $b=0$ and
the coordinate at $v_0(\lmb)$, one gets after simple transformations:
\vvn.6>
\begin{align*}
& \frac1{2\pi i}\;\frac{\Gm(\gm)}{\Gm(\al)\>\Gm(\bt)}
\,\int_{-\)i\)\infty\)-\eps}^{+\)i\)\infty\)-\)\eps}\!(-\)x)^s\;
\frac{\Gm(-\)s)\>\Gm(s+\al)\>\Gm(s+\bt)}{\Gm(s+\gm)}\ ds\,={}
\\[2pt]
&\ \ {}=\,(1-x)^{\gm-\al-\bt}\;\frac{\Gm(\gm)}{\Gm(\al)\>\Gm(\gm-\al)}\,\,
\int_{1\,}^{+\infty} t^\bt\>(t-1)^{\al-1}\>(t-x)^{\bt-\gm}\>dt\,={}
\\[2pt]
&\ \ \quad{}=\;\frac{\Gm(\gm)}{\Gm(\al)\>\Gm(\gm-\al)}
\,\,\int_{\!0}^{1\!}u^{\al-1}\>(1-u)^{\gm-\al-1}\>(1-ux)^{-\)\bt}\>du\,=\,
{}_2F_1(\al\),\bt\);\gm\);x)\,,
\cnn-.2>
\end{align*}
where
\vvn-.1>
$$
\al=-\)l_1/\ka\,,\qquad \bt=(z_1\]-z_2-l_1)/\ka\,,\qquad
\gm=(z_1\]-z_2-l_1\]+1)/\ka\,.
$$
Here it is assumed that \,${\Re\gm>\Re\al>0}$, \,${\Re\bt>0}$,
\,${0<\eps<\min\)(\)\Re\al\>,\Re\bt\))}$, and \,${-\pi<\arg\)(-\)x)<\pi}$,
\,${-\)\pi<\arg\)(1-x)<\pi}$. The second equality is obtained
by the change of integration variable $u\>=\>(t-1)/(t-x)$.
\vsk.2>
Theorems~\ref{dualhint} and \ref{dualqhint} exhibit the $\glkn$ duality for
\hint/s for $k=n=2$. The proofs of the theorems essentially involve explicit
formulae for Selberg-type integrals \(selb) and \(qselb)\). Those integrals
are associated with the Lie algebra $\slt$. To extend the duality of \hint/s
to the case of arbitary $k,n$ one needs to know suitable generalizations of
the Selberg integral associated with the Lie algebras $\slk$ for ${k>2}$.
For $k=3$ the required generalizations were obtained in \cite{TV5}\),
and similar ideas can be used to construct the required Selberg-type integrals
associated with the Lie algebras $\slk$ for ${k>3}$.

\section*{Acknoledgements}
The author was supported in part by RFFI grant 02\)\~\)01\~\)00085a \>and
CRDF grant RM1\~\)2334\)\~MO\)\~\)02.

\end{document}